\documentclass[11pt,a4]{article}

\usepackage{algorithm,algpseudocode}
\usepackage{amsmath,amsthm,amssymb}
\usepackage{geometry,subfig,hyperref,
epsfig,epstopdf,bbm, xcolor}

\newtheorem{theorem}{Theorem}[section]

\newtheorem{remark}[theorem]{Remark}

\numberwithin{equation}{section}
\newtheorem{exmp}[theorem]{Example}

\newcommand{\R}{\mathbb{R}}
\newcommand{\bxi}{\boldsymbol{\xi}}
\newcommand{\indep}{\rotatebox[origin=c]{90}{$\models$}}
\newcommand{\E}{\mathbb{E}}

\usepackage{enumerate}
\usepackage[shortlabels]{enumitem}

\newcommand{\balp}{\boldsymbol{\alpha}}

\title{A Dynamical Sparse Grid Collocation Method for Differential Equations Driven by White Noise}

\author{H. Cagan Ozen\thanks{Department of Applied Physics \& Applied Mathematics, Columbia University, New York, NY 10027, USA (\href{mailto:hco2104@columbia.edu}{hco2104@columbia.edu}, \href{gb2030@columbia.edu}{gb2030@columbia.edu}). }  \and Guillaume Bal\footnotemark[1] }

\date{}

\begin{document}
\maketitle

\begin{abstract}

We propose a sparse grid stochastic collocation method for long-time simulations of stochastic differential equations (SDEs) driven by white noise. The method uses pre-determined sparse quadrature rules for the forcing term and constructs evolving set of sparse quadrature rules for the solution variables in time. We carry out a restarting scheme to keep the dimension of random variables for the forcing term, therefore also the number of quadrature points, independent of time. At each restart, a sparse quadrature rule for the current solution variables is constructed based on the knowledge of moments and the previous quadrature rules via a minimization procedure. In this way, the method allows us to capture the long-time solutions accurately using small degrees of freedom. We apply the algorithm to low-dimensional nonlinear SDEs and  demonstrate its capability in long-time simulations numerically.

\end{abstract}

\pagestyle{myheadings}
\thispagestyle{plain}
\markboth{ H.C. Ozen and G. Bal}{DSGC}

\section{Introduction}
\label{sec:Introduction}

%
%
%
%
%
\medskip

Stochastic collocation methods have been popular tools for computing expectations of functionals of solutions to stochastic differential equations (SDEs) driven by white noise \cite{OB03, JK09, Kloeden, H01}. Undoubtedly, the most commonly used approach to compute expectations has been Monte Carlo (MC) method, which is a sampling technique in the random parameter space. However, random sampling-based methods are notorious for their slow convergence, which hinders their computational efficiency in simulations; especially in long-time simulations. Although several extensions have been proposed to speed up the convergence, their speed does not compete with fast convergence of spectral methods; see \cite{G08,H15,Kloeden,MT97,LeMK10} and references therein.    

In this work, we focus on stochastic collocation methods which use deterministic quadrature nodes in the random space to approximate expectations of functionals of solutions. These methods aim to achieve the ease of implementation of MC methods and fast convergence behavior of stochastic Galerkin methods at the same time. Similar to Galerkin methods, their convergence depends on the smoothness of the solution in the random input variables. Their effectiveness relies upon the dimensionality of the random parameter space and they work reasonably well if the stochastic system involves moderately large random dimensions. It has been shown, especially in uncertainty quantification literature, that they provide a strong alternative to MC methods for differential equations with time-independent random forcing; see e.g., \cite{XH05,Xiu07,BNT07,NTW08,NTW082, Xiu09, LeMK10}.

For equations driven by white noise, collocation methods with pre-determined quadrature rules have appeared in the recent literature \cite{GG98,Ger07,ZTRK14,ZTRK15, LV04, LL12}. Manuscripts \cite{ZTRK14, ZTRK15} combine a sparse grid collocation (SGC) method for white noise with weak-sense time integration methods to compute second order moments of the solutions of SDEs and stochastic partial differential equations (SPDEs). It has been proved and observed numerically in \cite{ZTRK14} that straightforward application of collocation methods leads a failure in long-time integration. The underlying reason of this failure is that in the presence of random forcing in time, the number of stochastic variables needed to maintain a prescribed accuracy increases with time. This means, for collocation-based methods, that the number of collocation points should grow with time. The manuscript \cite{ZTRK14} then introduces a recursive approach for long times on top of a straightforward implementation of SGC to compute moments of linear SPDEs; see also \cite{LMR97,MR04}. Similar long-time integration issues have been noted before in polynomial chaos expansion-based methods in \cite{WK06, GSVK10, BM13, OB16, OB17}. In contrast to collocation methods with pre-determined quadrature rules, optimal quantization methods \cite{PP05,LP06,PS11} aim to find optimal discrete approximations, e.g. Voronoi quantizations, to the solutions of SDEs, which are tailored for the underlying dynamics. These quantizations are usually obtained by approximating Brownian motion by a finite-dimensional random process, e.g. truncated Karhunen-Loeve expansion (KLE) of Brownian motion, and deriving systems of ordinary differential equations for the quantizers. Although resulting quadrature rules are adapted to the measures under consideration, slow convergence rates limit use of optimal quantization methods in numerical applications without further acceleration techniques \cite{PP05}.

A possible remedy to the long-time integration problem in the setting of SDEs is to exploit Markov property of the dynamics and allow the algorithm to forget about the past. This observation is crucially leveraged in our previous works \cite{OB16,OB17} to construct dynamical polynomial chaos expansions, called \textit{Dynamical generalized Polynomial Chaos} (DgPC), in terms of the solution variables and the random forcing; see also \cite{ZRTK2012} for a different approach. DgPC captures the solutions of SDEs and SPDEs accurately in the long-time with less computational time compared to MC methods; especially for complicated SPDE dynamics. In this work, using similar principles of \cite{OB16,OB17}, we propose a dynamical collocation-based method in time to alleviate long time integration of SDEs. We propagate optimal quadrature rules for the solution in time and use pre-determined quadrature rules for the random forcing. In this sense, the method can be considered as an extended combination of the proposed methods in \cite{ZTRK14,PP05}. The method uses a restarting scheme similar that of \cite{OB16,OB17} to construct sparse quadrature rules for the solution variables on-the-fly. It then estimates expectations of functionals of future solution variables by using sparse quadrature rules of the random forcing variables and the solution variables at the current time. By constructing quadrature rules with small number of nodes and employing frequent restarts, the algorithm can utilize small and time-independent degrees of freedom while maintaining accuracy in the long-time. We demonstrate the  efficiency of the proposed method numerically in several cases using low-dimensional nonlinear SDEs forced by white noise. 

The outline of the paper is as follows. Section \ref{sec:Quad} introduces the necessary background material on sparse quadrature methods. A simple stochastic collocation method is introduced in Section \ref{sec:SDE} and its drawbacks are discussed. Our methodology, called \textit{Dynamical Sparse Grid Collocation} (DSGC), is described in detail in Section \ref{sec:Formulation}. Numerical experiments for nonlinear SDEs are presented in Section \ref{sec:Numeric}. Some conclusions are offered in the last section.

\section{Sparse Quadrature Methods} \label{sec:Quad}

Approximation of an output $u(\bxi) \in \R^d$ of a model depending on a stochastic input variables $\bxi \in \R^K$ requires the solution of the corresponding forward propagation operator. This approximation procedure usually involves computation of expectations of the functional $u(\bxi)$. For instance, for polynomial chaos expansions (PCEs), the coefficients of the representation can be obtained by a spectral projection onto associated orthogonal polynomial basis $T_{\balp}(\bxi)$ and are given in the form of expectations $\E [u(\bxi) \, T_{\balp}(\bxi)]$ with respect to the probability measure of $\bxi$. The dimension $K$ of $\bxi$ is usually large in most applications, which in turn requires efficient estimations of high-dimensional integrals. In the following, we focus on quadrature-based collocation methods to compute expectations.

Assuming the probability distribution of the input variables $\bxi$ is known, e.g. it is given in the Askey family \cite{XK02}, and the components are independent, one can construct multi-dimensional quadrature rules by employing tensorization of one-dimensional quadrature rules. For the sake of brevity, we consider evaluation of $K$-dimensional integrals of the output $u(\bxi)$ with respect to a Gaussian measure $p_{\bxi}(\bxi) \propto \exp(-||\bxi||_2/2)$. 

In one dimension, we define the Gauss-Hermite quadrature rules $I_{Q}$ consisting of the weights and nodes $\{w^q, \xi^q \}_{q=1}^Q$, $Q \in \mathbb{N}$ and $\xi \in \R$: 
 \begin{align*}
 I_Q (u)(\xi) : = \sum_{q=1}^Q w^q \, u(\xi^q), 
 \end{align*}
where $\xi^q$'s are the roots of the Hermite polynomials $H_Q$ and the weights are given by $w^q = 1/(Q^2 (H_{Q-1}(\xi^q))^2)$. It is known that $I_Q$ is exact if $u$ is a polynomial of degree less than or equal to $2Q-1$; \cite{DR84}. 

For the multi-dimensional case, the integral $\E[u(\bxi)]$ can be approximated by the tensor product formula 
\[
\E[u(\bxi)] \approx I^{\otimes K}_Q  := \sum_{\alpha_1 =1}^{Q} \hdots \sum_{\alpha_K =1}^{Q} u(\xi_1^{\alpha_1}, \hdots,\xi_K^{\alpha_K}) \, w_1^{\alpha_1} \hdots w_K^{\alpha_K}.   
\]
Here we use the multi-index notation: $\balp = (\alpha_1, \hdots, \alpha_K) \in \mathbb{N}^{K}$ with $|\balp| = \sum_{k=1}^K \alpha_k$. We denote by $Q_{\bxi}$ the total number of resulting quadrature nodes. Approximations based on this tensor product suffer from curse of dimensionality and computational costs scale exponentially with dimension $K$, i.e. $Q_{\bxi}= Q^K$.

If the dimension of the random variables is moderately high, a sparse quadrature rule first proposed by Smolyak, can be used to reduce the number of quadrature nodes while maintaining the accuracy \cite{smolyak}. Following \cite{WW95,GG98}, we write the sparse grid approximation to the multi-dimensional integral with the level $\lambda$ 
\begin{align} \label{eq:Smolyak}  
\E[u] \approx \sum_{ \lambda \leq |\balp| \leq \lambda + K -1} (-1)^{\lambda+K-|\balp|-1} {K-1 \choose |\balp| -\lambda}  ( I_{\alpha_1} \otimes \hdots \otimes I_{\alpha_K} )(u),
\end{align}
where $\balp \geq \mathbf{1}$. This quadrature rule is exact for multivariate polynomials of total degree up to $2\lambda-1$ and greatly reduces the number of evaluations compared to the tensor product rule above \cite{GG98,NR99}.  In this work, we employ isotropic Smolyak sparse grid quadrature rules for Gaussian measures, meaning that the level $\lambda$ is the same for each dimension. We also note that the weights of this sparse quadrature rule can be negative. 

For stochastic equations with random input data, sparse grid collocation methods have been shown to be effective, especially equations with time-independent forcing and parameters, and in moderately high dimensions; see \cite{XH05,Xiu07,Xiu09,BNT07,NTW08,NTW082} and references therein. For stochastic equations with time-dependent noise, there are mainly two methods which use pre-determined set of quadrature nodes and weights: cubatures on Wiener space \cite{LV04, LL12} and sparse grid Smolyak quadrature rules \cite{GG98,Ger07,ZTRK14,ZTRK15}. In \cite{ZTRK14}, applications of sparse grid collocation methods to equations driven by white noise have been studied in conjunction with Monte-Carlo methods and it is noted that sparse grid collocation works only for short time integration. Motivated by this observation, in the next section, we introduce a sparse, deterministic collocation method based on quadrature rules for SDEs driven by white noise using a spectral expansion of the white noise and confirm the findings in \cite{ZTRK14}.

\section{A stochastic collocation method for SDEs}  \label{sec:SDE}
 We consider the  following $d$-dimensional SDE:
\begin{align} \label{eq:sdeD}
d u(t)= \mathcal{L}(u(t)) \, dt + \sigma(u(t)) \, dW(t), \quad u(0)=u_0, \quad t \in [0,T],
\end{align}
where $\mathcal{L}(u)$ is a general function of $u \in \R^d$ (and possibly other deterministic or stochastic parameters), $W(t)$ is a Brownian motion, $\sigma(u)$ is the amplitude which depends on the solution, and $u_0$ is an initial condition with a prescribed probability density. The stochastic integral is considered in the Ito sense. We assume that the solution exists and is a second order stochastic process. We also use the notation $W(t,\omega)$ for Brownian motion $W(t)$.

We introduce a sparse quadrature-based collocation method for \eqref{eq:sdeD} as follows. First, we consider a time discretization of the interval $[0,T]$
\[
 \tau_i = i \, \Delta \tau, \quad i=0,\dots,M_{T},
\]
where $\Delta \tau = T/M_{T}$. Then we approximate the solution $u(t; u_0, \{W({\tau}); \tau \leq t\})$ of \eqref{eq:sdeD} by the Euler method \cite{Kloeden}
\begin{align*}
u(\tau_{i+1}) & = u(\tau_{i}) + \mathcal{L}(u(\tau_{i})) \, \Delta \tau + \sigma(u({\tau_i})) \, (W(\tau_{i+1}) - W(\tau_{i})).
\end{align*}

Given a complete countable orthonormal system $m_k(t) \in L^2[0,T]$, $k \in \mathbb{N}$, we project the Brownian motion $W(t,\omega)$ onto $L^2 [0,T]$ by defining $\xi_k(\omega) := \int_0^T m_k(t)\, dW(t,\omega) $. Then, the random vector $(\xi_1,\xi_2,\ldots)$ consists of a countable number of  independent and identically distributed standard Gaussian random variables, and the convergence 
\begin{align} \label{eq:conv_W}
 \E \left[ W(t) - \sum_{k=1}^K  \xi_k \, \int_0^t m_k(s) \, ds \right]^2  \rightarrow 0 , \quad K \rightarrow \infty,
\end{align}
holds for all $t \leq T$; see~\cite{HLRZ06, Luo,MNGK04}.
This convergence property allows us to approximate the increments $(W(\tau_{i+1}) - W(\tau_{i}))$ of the random forcing by the finite dimensional random variable $\boldsymbol{\xi}= (\xi_1,\xi_2,\hdots,\xi_K)$:
\begin{align} \label{eq:sde_integral}
{u}(\tau_{i+1}; u_0,\bxi) & = {u}(\tau_i; u_0,\bxi)  + \mathcal{L}({u}(\tau_{i}; u_0,\bxi)) \, \Delta \tau + \sigma(u(\tau_i;u_0,\bxi)) \sum_{k=1}^K \xi_k \int_{\tau_{i}}^{\tau_{i+1}} m_k(s) \, ds,
\end{align}
where $ {u}(\tau_0; u_0,\bxi)  =u_0$. 

Now, let $\{w_{0}^p, u_0^p\}_{p=1}^{Q_{u_0}}$ and $\{\mathbf{w}^q, \bxi^q\}_{q=1}^{Q_{\bxi}}$ be pre-determined quadrature rules for the random variables $u_0$ and $\bxi$. For instance, we can use Smolyak sparse grid quadrature \eqref{eq:Smolyak} for Gaussian $\bxi$ and for $u_0$, any accurate quadrature rule can be considered. Here these quadrature rules denote any enumerations of their multi-dimensional versions. Then equation \eqref{eq:sde_integral} naturally gives rise to a non-intrusive approximation to \eqref{eq:sdeD} by the following equation: 
\begin{align} \label{eq:particleSys}
{u}(\tau_{i+1}; u_0^p,\bxi^q) & = {u}(\tau_i; u_0^p,\bxi^q)  + \mathcal{L}({u}(\tau_{i}; u_0^p,\bxi^q)) \, \Delta \tau + \sigma(u(\tau_i;u_0^p,\bxi^q)) \sum_{k=1}^K \xi^q_k \int_{\tau_{i}}^{\tau_{i+1}} m_k(s) \, ds,
\end{align}
for $p=1,\hdots,Q_{u_0}$ and $q=1,\hdots,Q_{\bxi} $. Equation \eqref{eq:particleSys} dictates the evolution of the initial particles $u_0^p$ under the trajectories of the forcing particles $\bxi^q$. In contrast to Monte Carlo methods, the random forcing is sampled deterministically and approximated by its finite-dimensional approximation via the spectral projection \eqref{eq:conv_W}, and the samples of $u_0$ are taken as quadrature points. Thus, the method is sample-error free.  

\begin{remark} \rm 
The method introduces three level of approximations. First, SDE is discretized in time. Then, Brownian motion increments are approximated by their finite-dimensional approximations. Finally, the semi-discrete equation \eqref{eq:sde_integral} is approximated by its fully discrete version \eqref{eq:particleSys} using quadrature rules. Hence, there are three degrees of freedom that are of interest: $\Delta \tau$, $K$, and $\lambda$. Note also that although we used the Euler method in the formulation, this is not required. Any higher order method can be used to discretize the SDE in time.  
\end{remark}

\begin{remark} \rm For any fixed $K$, the finite dimensional approximation \eqref{eq:conv_W} of Brownian motion entails a continuous finite-variation process on $[0,T]$. Thus, the integral with respect to this finite-variation process can be understood in the Stieltjes sense. Then the main questions are when and in what sense the approximate solution \eqref{eq:particleSys} converges to the true solution. Here we are only intereseted in the numerical efficiency of the method and we refer to \cite{PP05,LP06,PS11} for theoretical discussions on the convergence of approximations for SDEs with smooth coefficients in a similar setting. Nevertheless, our numerical experiments show clear convergence in moments; see Section \ref{sec:Numeric}. 
\end{remark}

The approximate solution $u(t;u_0^p,\bxi^q)$ of \eqref{eq:particleSys} readily entails approximations to statistical moments by computing 
\[
\E[u(t; u_0,\bxi)^{\balp}] \approx \sum_{p=1}^{Q_{u_0}} \sum_{q=1}^{Q_{\bxi}} w_0^p \, \mathbf{w}^q \, u(t; u_0^p, \bxi^q)^{\balp}, 
\]
where $\balp $ is a multi-index, see \eqref{eq:multiindex}.

A similar collocation strategy in a Monte-Carlo setting using weak-integration is employed in \cite{ZTRK14,ZTRK15} to compute first two moments of the solution. It is noted that the efficiency of the collocation strategy depends on the strength of the noise and the length of the time interval. Indeed, in order to maintain a prescribed accuracy the expansion \eqref{eq:conv_W} requires the number of stochastic variables to be increasing with time. Thus, the number of quadrature points needed to maintain an accuracy becomes quickly overwhelming for long times.

\begin{figure}[!htb]
\centering
\subfloat[Relative error versus $K$]{
\includegraphics[scale=0.4]{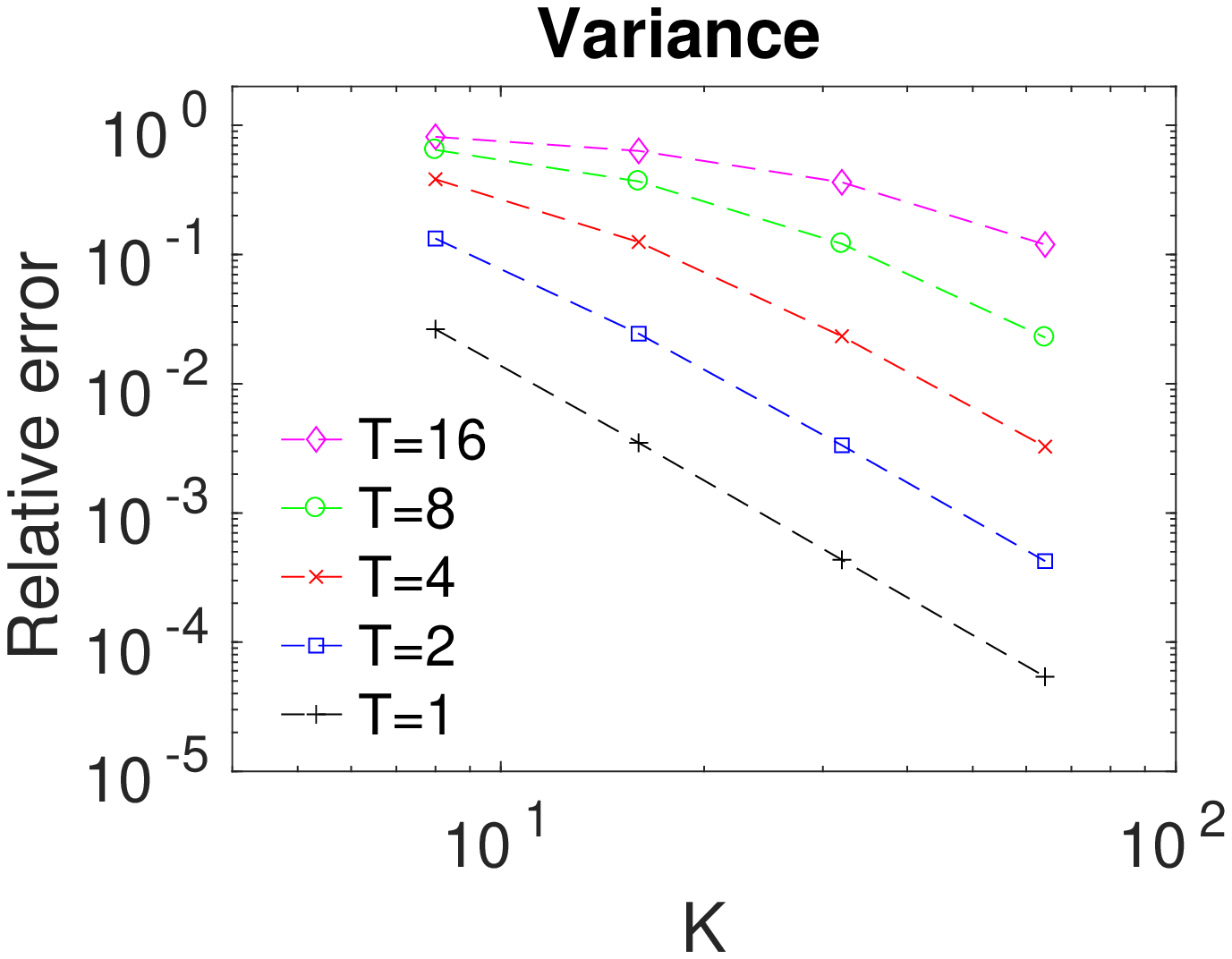}} 
\subfloat[Relative error versus $T$]{
\includegraphics[scale=0.4]{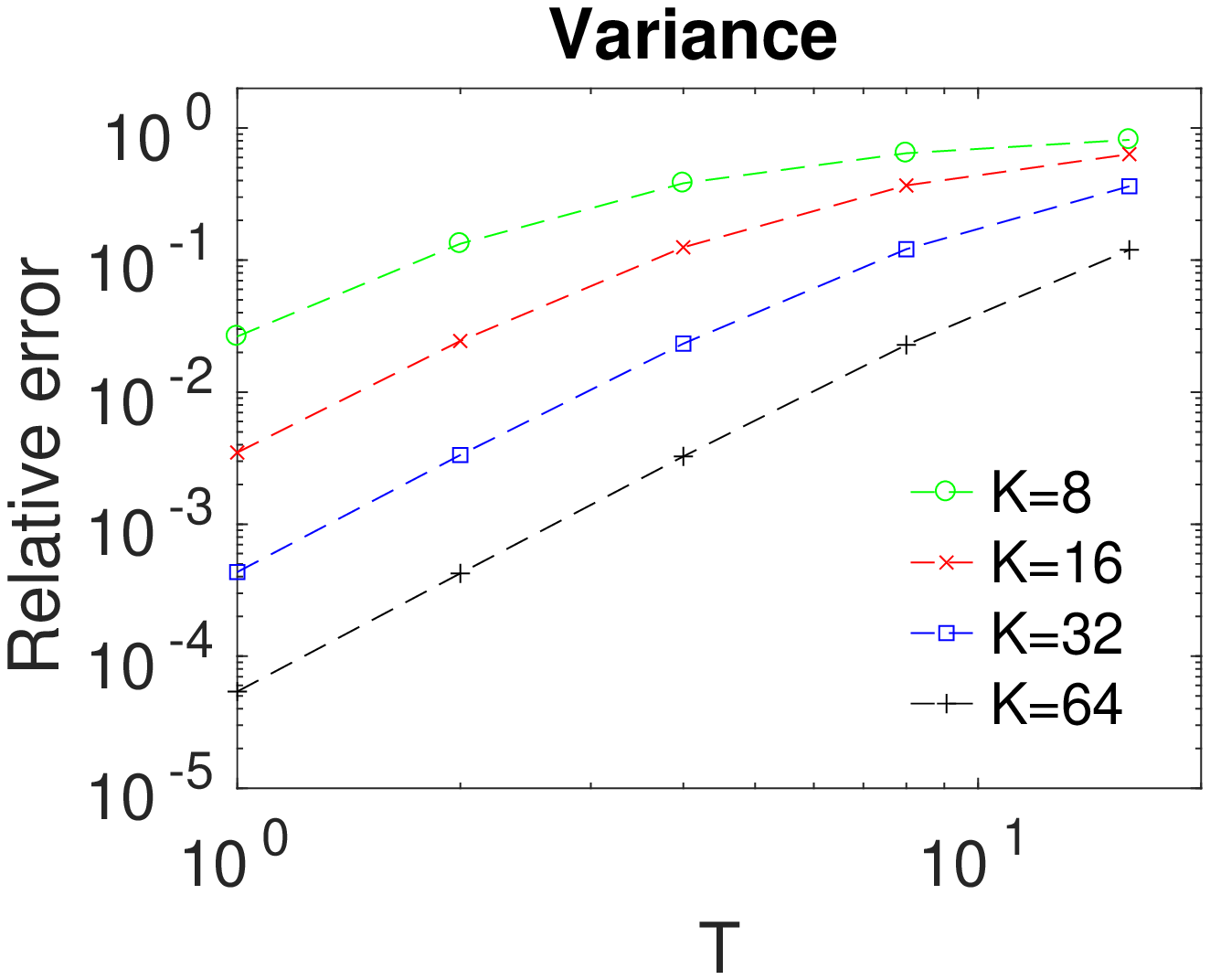}} 
\caption{Relative errors of the variance for different times $T=1,2,4,8$ and $16$, and different number of random variables $K=8,16,32$ and $64$ for the Ornstein-Uhlenbeck process.}
\label{fig:OU}
\end{figure}

Here is a simple motivating demonstration of the long-time integration problem in case of the Ornstein-Uhlenbeck process. We set $\mathcal{L}(u)= 10(0.1-u)$ and $\sigma=4$, and take a deterministic initial condition $u_0=1$. We also take $K=8, 16, 32, 64$, and consider different final times $T=1,2,4,8,16$. Since the solution stays Gaussian, we use a sparse Gauss-Hermite rule for $\bxi$ with level $\lambda_{\bxi} = 1$. To make a fair comparison, we use the same time discretization method (second order Runge-Kutta method) with the time step $ \Delta \tau= 1$E-$3$ in each scenario. Figure \ref{fig:OU} shows that as the time increases from $T=1$ to $T=16$, the convergence behavior in the variance of this simple method decreases from $O(K^{-3})$ to $O(K^{-1})$. This clearly indicates that the degrees of freedom required for this simple collocation method to maintain a desired accuracy should increase with time.

\section{Proposed methodology}  \label{sec:Formulation}

\subsection{Motivation} 

As our previous works \cite{OB16,OB17} and other manuscripts \cite{MR04,HLRZ06,ZRTK2012,ZTRK14,ZTRK15} noted before,  the main difficulty for differential equations driven by white-noise forcing that should be addressed is that the number of variables need to be accounted for grows in time. In the settings that dynamics satisfy a Markov property, any adapted algorithm can and should utilize this fact to ``forget" the past variables and as a consequence keep the dimension of the random variables independent of time. Basically, this formulation leverages the intrinsic sparsity of the dynamics in the sense that a sparse representation of future times can be obtained using the solution at the current time and the corresponding random forcing variables. In \cite{OB16,OB17}, utilizing these observations, we introduced a PCE-based method which utilizes a restarting procedure to compute orthogonal polynomials of the solution dynamically on-the-fly and constructs evolving polynomial expansions to represent future solution in a reasonably sparse representation. The method performed on par with standard MC methods in most cases and, in case of complex dynamics, it can be executed faster than standard MC \cite{OB17}. Although, PCE-based restart methods provide a viable means to keep the number of random variables in reasonable levels and to compute statistical properties of long-time evolutions for fairly complicated equations, the methods in \cite{OB16,OB17} are intrusive. In this connection, we present a non-intrusive method which exploits the same forgetting strategies to construct dynamical quadrature-based collocation rules for the solution in time. For other forgetting strategies, we refer the reader to \cite{MR04,LL12,ZRTK2012,ZTRK14,ZTRK15}.

\subsection{Formulation} 

The solution $u(t,\omega)$ of the $d$-dimensional SDE \eqref{eq:sdeD} is a nonlinear functional of the initial condition and countable number of variables $\bxi$. The restart methods in \cite{OB16,OB17} crucially depend on the following two observations: (i) for sufficiently small time lags $\varepsilon>0$, the solution $u(t+\varepsilon)$, can be efficiently captured by low-order polynomials in $u(t)$; and (ii) the solution $u(t+\varepsilon)$ depends on Brownian motion only on the interval $[t,t+\varepsilon]$ by Markov property. Based on these observations, the algorithms in\cite{OB16,OB17} introduced an increasing sequence of restart
times $0 < t_j < t_{j+1} < T$, with $\Delta t= t_{j+1} - t_j$, and constructed a new PCE basis at each $t_j$ based on the solution $u(t_j )$ and all additional random variables $\bxi$. 

In the following, we propose to utilize quadrature rules of the solution $u(t)$ and the variables $\bxi$ to represent the solution $u(t+\epsilon)$ at future times with a sufficiently small number of quadrature points. By the fact (i), $u(t+\epsilon)$ can be captured by low order polynomials in $u(t)$, therefore, the quadrature level required to integrate polynomials in $u(t)$ can be selected small. Moreover, the number of quadrature points for $\bxi$ on each time interval $[t_j,t_{j+1}]$ can also be made small by selecting a short time horizon. Notice also that by the properties of Brownian motion, the quadrature rule for $\bxi$ can be read from tables or computed only once in the offline stage. The challenge is then to compute efficient, sparse quadrature rules for the solution $u(t)$ in time. These rules are not straightforward to compute and have to be computed online since the probability distribution of $u(t_j)$ is arbitrary and evolving.  

Let a sequence of increasing restart times $t_j$ be given.  We denote by $u_j$ the approximation of $u(t_j)$ given by the algorithm at $t_j$ and by $\bxi_j=(\xi_{j,1},\hdots,\xi_{j,K})$ the variables for the random forcing on the interval $[t_j,t_{j+1}]$. Suppose for now that the sparse quadrature rules $\{w^p_j, u_j^p \}_{p=1}^{Q_{u_j}}$ and $\{\mathbf{w}^q_j, \bxi_j^q\}_{q=1}^{Q_{\bxi_j}}$ have already been constructed for $u_j$ and $\bxi_j$ at time $t_j$, respectively. An analog of the equation \eqref{eq:particleSys} can be written for the approximate solution $u_{j+1}(t; u_j,\bxi_j)$ as
\begin{align}  \label{eq:upq}
\nonumber
{u_{j+1}}(\tau_{j,i+1}; u_j^p,\bxi_j^q) & = {u}(\tau_{j,i}; u_j^p,\bxi_j^q)  + \mathcal{L}({u_{j+1}}(\tau_{j,i}; u_j^p,\bxi_j^q)) \, \Delta \tau \\ 
& \quad + \sigma(u_{j+1}(\tau_{j,i};u_j^p,\bxi_j^q)) \sum_{k=1}^K \xi^q_{j,k} \int_{\tau_{j,i}}^{\tau_{j,i+1}} m_{j,k}(s) ds,
\end{align}
where $m_{j,\cdot}(t)$ is a complete orthonormal system for $L^2[t_j,t_{j+1}]$ and ${\tau_{j,\cdot}}$'s denote a time discretization for the interval $[t_{j},t_{j+1}]$. 

The question here is how to construct an efficient quadrature rule for the next solution variable $u_{j+1}$ using \eqref{eq:upq}. The evolution of the particles $\{u_j^p, \bxi_j^q\}$ via the equation \eqref{eq:upq} entails a set $\{u^{p \times q}_{j+1}\}$ of particles of the approximate solution $u_{j+1}$, i.e quadrature nodes at $t_j$ follow the trajectories of the dynamics and give rise to an initial set of nodes at $t_{j+1}$. The challenge is to find a small subset of these nodes and corresponding weights which accurately integrates polynomials in $u_{j+1}$. Following \cite{AGPRH12,AGPRH14} and \cite{EB2015}, we construct such a sparse quadrature for $u_{j+1}$  by using the following $L^1$-minimization procedure. 
 
We define multi-indices $\balp$ up to degree $N$ 
\begin{align} \label{eq:multiindex}
J_{d,N} := \{ \balp=(\alpha_1,\ldots,\alpha_d) \in \mathbb{N}_0^d, \, |\balp|\leq N\}.
\end{align} 
and let $M:= |J_{d,N}| = {d+N \choose{ d}}$ be the total number of elements in the set. 
Let the particles $u^{p \times q}_{j+1}:= u_{j+1}(t_{j+1}; u_j^p, \bxi_j^q)$, where $p=1,\hdots, Q_{u_j}$ and $q=1,\hdots, Q_{\bxi_j}$, serve as an initial set of quadrature nodes for $u_{j+1}$. Let also the set
$\{ T_{\balp}(u) \, : \balp \in \mathcal{J}_{d,N} \}$
be any orthonormal basis of polynomials up to degree $N$ in dimension $d$. Then to extract a sparse quadrature rule, following \cite{AGPRH12,AGPRH14}, we solve the convex optimization problem: 
\begin{align} \label{eq:optL1}
\min_{w_{j+1} \in \R^{\tilde{Q}_{u_{j+1}} } } || w||_1 , \quad \mbox{subject to  } A w  = b, 
\end{align}
where $w \in \R^{\tilde{Q}_{u_{j+1}}}$ with $\tilde{Q}_{u_{j+1}} = Q_{u_j} \times Q_{\bxi_j}$, and the constraints $A w = b $ necessitate the exactness of the quadrature rule up to degree $N$. We enumerate the basis $ T_{\balp}(u)$ and denote it by $\{T_{k}(u) \, : \, k =0,\hdots , M\}$, then define
\begin{align*}
A := \begin{bmatrix}
T_{0}(u_{j+1}^1)  & \ldots  & T_{0}(u_{j+1}^{\tilde{Q}_{u_{j+1}}}) \\ 
\vdots &     & \vdots\\ 
T_{M}(u_{j+1}^1) & \ldots &  T_{M}(u_{j+1}^{\tilde{Q}_{u_{j+1}}}) \\ 
\end{bmatrix},
\end{align*}
and the right-hand side vector consisting of the exact moments 
\[
 b:= \begin{bmatrix}
 \E [T_{0}(u_{j+1})] & \ldots & \E [T_{M}(u_{j+1})] \end{bmatrix}^T. 
\]
Here we assume that the moments $\E [u_{j+1}^{\balp}]$, for each $t_{j+1}$ and $|\balp| \leq N$, are finite. Also, we typically have that $M$ is much smaller than $\tilde{Q}_{u_{j+1}}$. Then a sparse subset denoted by $\{w^{p}_{j+1}, u^{p}_{j+1}\}_{p=1}^{Q_{u_{j+1}}}$ with $Q_{u_{j+1}} \leq M \ll \tilde{Q}_{j+1}$, can be extracted having at most $M$ nodes from the solution of the optimization procedure \eqref{eq:optL1}; see section \ref{sec:Imp}. Once a quadrature rule for $u_{j+1}$ is constructed, the algorithm restarts and evolves the nodes on the next time interval according to \eqref{eq:upq}.

\begin{remark} \rm
For probability measures on $\R^d$, the existence of an exact quadrature rule with positive weights is guaranteed; \cite{Tchakaloff, Putinar,D67,gautschi}. In general, we do not enforce positivity condition for the weights $w$ since a sparse optimal solution with positive weights may not exist; see Section \ref{sec:Imp} for further details. We note that \eqref{eq:optL1} is not the only construction to find optimal quadrature rules; see also \cite{GW69, D67, AGPRH12,AGPRH14} and references therein. Furthermore, the convergence of exact quadrature rules for compactly supported probability measures has been studied extensively and the results can be found in classical literature \cite{gautschi,DR84}.
\end{remark}

\begin{remark} \rm
We do not tensorize quadrature rules for each component of $d$-dimensional random vector $u_{j}$. It is quite possible that components of $u_j$ exhibit correlation; therefore the tensorization is not optimal. However, since we impose constraints on multivariate moments $\E[u_{j}^{\balp}]$, the algorithm automatically establishes a quadrature rule for the full vector $u_j$. Moreover, if the dimension $d$ is high, the number of constraints $M$ in \eqref{eq:optL1} can be reduced by considering a sparse version of the multi-index set $\mathcal{J}_{d,N}$; see \cite{OB16, OB17, HLRZ06, BS11}. 
\end{remark}

\begin{remark} \rm  This remark concerns the differences and similarities of our approach to the methods in \cite{ZTRK14,ZTRK15} and optimal quantization methods in \cite{PP05,LP06}. 

Our method uses pre-determined quadrature rules for the random forcing variables and does not a have sampling error similar to the method in \cite{ZTRK14}. The main difference is that the paper first discretizes the stochastic equations in time and then considers quadrature rules for the random forcing variables in each time-step, i.e. the dimension of randomness depends on the resolution of the time-integration and might grow rapidly with fine discretization. In contrast, our method discretizes in the random space by considering the $L^2$-approximation \eqref{eq:conv_W} of Brownian motion with a fixed degree of freedom $K$. Although a different restarting mechanism is used in \cite{ZTRK14}, their method can only compute moments up to second order, whereas approximations to higher order moments are available in our method by quadrature rules provided approximations to higher moments converge. 

Our method finds optimal quadrature rules adapted to the evolving solution in a similar sense to optimal quantization methods. The quantization methods in \cite{PP05,LP06} aim to discretize the paths of an infinite dimensional randomness by a random vector in finite dimension. Finite dimensional approximations are obtained by the solution of a minimization procedure and are typically given by Voronoi cell collocation. For Brownian motion, quantizers based on its KLE are considered. Then the evolution of these quantizers for SDEs is obtained by solving a corresponding integral equation, which is similar to \eqref{eq:upq}. The convergence order is logarithmic in the number of quantizers, which is a poor rate of convergence for practical applications. In contrast, our method utilizes Gaussian-type quadratures tailored for the solution and the random forcing, which leverage the smoothness of the response to provide fast convergence. Moreover, frequent restarts allow us to mitigate dimensionality and use low-dimensional approximations to  Brownian motion forcing. 
\end{remark}

Our methodology is given in Algorithm 1 below.   
\begin{algorithm}
\caption{Dynamical Sparse Grid Collocation (DSGC) method for SDEs}
\label{alg:DSGC}
\begin{algorithmic}
\\ Decompose the time domain $[0,T] = [0,t_1] \cup \ldots \cup [t_{n-1},T]$
\\ Select a time-integration method
\\ Initialize the degrees of freedom $K,N$
\\ Compute quadrature rules for $\bxi_0$ and $u_0$
\For{each time-step $t_j$, $j>0$,}
\State evolve the quadrature nodes $u_{j-1}^p$ and $\bxi_{j-1}^q$ by \eqref{eq:upq}
\State obtain the nodes $u^{p \times q}_j$,  $p=1,\hdots, Q_{u_{j-1}}$ and $q=1,\hdots, Q_{\bxi_{j-1}}$
\State solve the optimization procedure \eqref{eq:optL1}
\State extract a sparse quadrature rule $\{w_{j}^p,u_{j}^p \}_{p=1}^{Q_{u_j}}$
 \EndFor 
\end{algorithmic}
\end{algorithm}

\subsection{Implementation} \label{sec:Imp}

\subsubsection{Offline stage}

In the offline stage, certain quadratures need to be computed. First, we compute sparse quadrature rule for Gaussian $\bxi_0$ by using the Smolyak sparse grid with the level $\lambda_{\bxi_0}$, which builds upon the standard 1D Gauss-Hermite rule; see Section \ref{sec:Quad}. By the independent increment property of Brownian motion, all $\bxi_j$, $j \geq 0$, may have the same quadrature rule assuming $K$ is fixed. Note that although it is not necessary, we keep the number of variables in $\bxi$ the same throughout the evolution. If the distribution of the initial condition $u_0$ is known, a sparse Gauss quadrature is constructed with the level $\lambda_{u_0}$. If its distribution is arbitrary, then the optimization procedure \eqref{eq:optL1} can be used with Monte Carlo initialization. We make use of the C++ library ``UQ Toolkit" to compute Gauss rules \cite{DNPKGL04}.

\subsubsection{Moments and orthogonal polynomials}
\label{sec:orth}

At the restart time $t_{j+1}$, estimation of the right-hand side vector in the constraints in \eqref{eq:optL1} requires the calculation of multivariate moments $\E[u_{j+1}^{\balp}]$, where $|\balp| \leq N$ and $u_{j+1}= (u_{j+1,1}, \hdots, u_{j+1,d})$. As noted before in Section \ref{sec:SDE}, these moments are computed using the available quadrature rules $\{w_{j}^p , u_{j}^p\}_{p=1}^{Q_{u_j}}$ and $\{\mathbf{w}^q_j, \bxi_j^q\}_{q=1}^{Q_{\bxi_{j}}}$ from time $t_j$:
\begin{align*}
\E[u_{j+1}^{\balp}] = \E_{(u_{j},\bxi_j)}[u_{j+1}^{\balp }(u_j, \bxi_j)] \approx \sum_{p=1}^{Q_{u_j}} \sum_{q=1}^{Q_{\bxi_j}} w_j^p \, \mathbf{w}_j^q \, \prod_{i=1}^d (u_{j+1,i}(u_j^p,\bxi_j^q))^{\alpha_i}.
\end{align*}

The optimization procedure does not depend on the particular choice of the polynomials $\{T_{\balp}(u) \, : \, \balp \in \mathcal{J}_{d,N} \}$, e.g. even monomials can be used. However, the choice of $T_{\balp}$ certainly affects the condition number of the constraint matrix, which in turn affects the stability of the numerical minimization algorithm. To better condition the constraints and improve convergence of the optimization algorithms, we make use of couple of linear transformations as preconditioning steps. Similar transformation techniques are applied  before in moment-constrained maximum entropy problem \cite{A10}. 

A linear transformation is applied to the solution $u_j$ so that it becomes mean zero. Then a further transformation makes its components uncorrelated, i.e. its covariance matrix becomes identity. Even with these transformations, a scaling issue related to the moments arises in the constraint equations. For instance, for a standard Gaussian random variable $\xi$, we have $\E[\xi^{12}]/\E[\xi^2] = 10395$, i.e. the twelfth moment is larger than the second moment by $5$ orders magnitude. To alleviate this scaling issue, we further scale $u$ by its largest moment so that the maximum moment becomes $1$. 

A more direct preconditioning can be applied by a judicious selection of the orthonormal basis $T_{\balp}$. A basis can be selected using an educated guess in the offline stage, which does not require any online computation. However, since the measure of $u_j$ is evolving in time, this may not be optimal in the long-time in Algorithm 1. An optimal choice for $T_{\balp}$ is the set of polynomials which are orthogonal with respect to the measure of $u_j$. Unfortunately, a corresponding orthogonal system for $u_j$ is not available a priori in our algorithm, but it can be computed online if a further preconditioning is required. Although the computation of orthogonal polynomials is an ill-posed problem, a Gram-Schmidt procedure based on the knowledge of multivariate moments can be used in the computation. In numerical simulations, the choice of the orthonormal system will be explicitly stated. We refer to \cite{gautschi,OB16,OB17} and references therein for detailed discussions on how to compute orthogonal polynomials of $u_j$ dynamically. 


Here is a simple demonstration of the effects of these transformations. Let $\xi_1$ and $\xi_2$ be two independent $N(3,1)$ variables and the maximum degree be $N=8$. Then, the number of constraints becomes $M=45$. We use $500$ samples from normal distribution to initialize the optimization procedure and keep the samples same for each scenario to make a fair comparison. A sparse quadrature rule with $M$ nodes is extracted according to the algorithm discussed in the next section, and, afterwards, the right-hand side vector $b$ is computed numerically using this quadrature rule to check the accuracy. The numerical approximation of $b$ is denoted by $\tilde{b}$ in the following. 

Table \ref{table:Scaling} shows condition numbers of the linear system in \eqref{eq:optL1} and the accuracy  $||b-\tilde{b}||_{\infty}/ ||b||_{\infty}$ of the associated quadrature rule. First two scenarios use monomials as the polynomial basis $T_{\balp}$ and the last one uses Hermite polynomials, which are the associated orthogonal system in this example. Note finally that condition numbers are independent of the sparse extraction procedure. Clearly, scaling transformations or a careful selection of the polynomials basis leads to at least $5$-digit gain of accuracy in this example.  

\begin{table}[!htb]
\begin{minipage}[b]{1\linewidth}
\centering
\renewcommand{\tabcolsep}{0.1cm}
\renewcommand{\arraystretch}{1.2}
 \begin{tabular} { |r | r |r|r|}
  \hline 
   & Without scaling & With scaling  On &  Hermite poly.\\ \hline
   cond(A) & 7.01E+9 &  1.91E+3 & 2.98E+2\\ 
   $||b-\tilde{b}||_{\infty}/ ||b||_{\infty}$ & 3.25E-8 &  3.95E-13 & $3.19$E-$14$  \\
    \hline  
  \end{tabular}
  \caption{The accuracy of quadrature rules for two independent Gaussian variables using different transformations.}
  \label{table:Scaling}
\end{minipage}
\end{table}

\subsubsection{Sparse quadrature rule}

Algorithm 1 constructs dynamical quadrature rules in time for the solution $u_j$. However, the implementation of the optimization algorithm \eqref{eq:optL1} to construct an efficient quadrature rule is not straightforward. From \eqref{eq:optL1}, we observe that at each restart $t_j$, the optimization procedure is initialized with $\tilde{Q}_{u_{j}} = Q_{u_{j-1}} \times Q_{\bxi_{j-1}}$ number of nodes for $u_j$. Therefore, the number of quadrature nodes may grow with the number of restarts provided the solution of the optimization procedure is not already sparse. The challenge is then to compute a sparse quadrature rule containing a smaller set of nodes and weights while keeping the exactness of the original quadrature rule. To this end, after finding the optimal solution of \eqref{eq:optL1}, we further employ an extraction routine. 
 
A straightforward sparsification of the optimal solution of \eqref{eq:optL1} would be cutting-off the weights that are greater than a certain threshold. Depending on the numerical minimization algorithm, this may not be possible. As discussed in \cite{AGPRH12,AGPRH14}, an application of simplex algorithm yields a sparse optimal solution whereas interior-point methods give fully populated solution \cite{N06}. This behavior is due to the fact that the simplex algorithm acts on the vertices of a polytope whereas interior point methods result in sparse vectors only in the limit, which is not achieved in numerics. We choose to use CVX, a package for specifying and solving convex programs \cite{cvx, gb08}. Under the hood, CVX uses SDPT3 which employs interior-point methods to compute the optimal solution \cite{SDPT3}. The following procedure is used to extract a sparse quadrature rule; see \cite{AGPRH12,AGPRH14,N06}. 

At time $t_j$, the constraints matrix $A$ is of dimension $M \times \tilde{Q}_{u_j}$, where $\tilde{Q}_{u_j}$ is much bigger than $M$. Then, the dimension of the nullspace of $A$ is at least $\tilde{Q}_{u_j} - M$. The key observation is that any vector $z \in \mbox{null}(A)$ can be added to the weights vector without changing the equality constraints, i.e. $A(w+z) = b$. Thus, by selecting vectors carefully, we can construct an iterative routine to make most of the weights zero. Based on these observations, we follow the approach given in \cite{AGPRH12,AGPRH14,N06} and employ the following routine, Algorithm \ref{alg:sp_ext}, at each restart $t_j$.  

\begin{algorithm}
\caption{Sparse Quadrature Extraction Routine for \eqref{eq:optL1}}
\label{alg:sp_ext}
\begin{algorithmic}
\\ Initialize with the optimal weights $w \in \R^{\tilde{Q}}$ 
\Repeat 
\State find the indices $\mathcal{N}:= \{ k \in \{1,\hdots, \tilde{Q} \} \, : \, w_k=0 \}$
\State find $z \in \mbox{null}(A)$ with $z_k=0$, $k \in \mathcal{N}$
\State set $\beta = \min \{ |\frac{w_k}{z_k}| \, : k \not \in \mathcal{N}, \, \mbox{sign}(z_k) \neq  \mbox{sign}(w_k)\}$
\State set $w = w + \beta z$
 \Until  the number of nonzeros in $w$ is less than or equal to $M$
\end{algorithmic}
\end{algorithm}

Algorithm \ref{alg:sp_ext} allows us to find a sparse quadrature rule  $\{w_j^p, u_j^p \}_{p=1}^{Q_{u_j}}$ with the number of nodes $Q_{u_j}$ satisfying $Q_{u_j} \leq M \ll \tilde{Q}_{u_j}$. Thus, the number of nodes can be made independent of time and frequent restarts can be used in Algorithm \ref{alg:DSGC}. A one way to find a vector $z$ in the nullspace of $A$ is to compute a basis for the nullspace. In the implementation, we make use of the QR method to quickly select a vector at each iteration. Finally, we note that this routine does not necessarily yield nonnegative weights. 

An application of this procedure to the two dimensional Gaussian random variable discussed before in Table \ref{table:Scaling} reduces the size of the quadrature rule from $500$ to $45$ while the accuracy remains almost the same as $3.19$E-$14$.

\section{Numerical Experiments} \label{sec:Numeric}

In this section, we present several numerical simulations of our algorithm using low-dimensional nonlinear SDE models. We refer to \cite{BM13,OB16,OB17} for details on the following numerical examples. 

For the rest of the section, $T \in \R$ stands for the endpoint of the time interval while $\Delta t= T/n$ denotes the time-step after which restarts occur at $t_j= j \Delta t $. Furthermore, following~\cite{Luo,HLRZ06,OB16,OB17}, we choose the following orthonormal basis for $L^2[t_{j},t_{j+1}]$:
\begin{align*}
m_{j,1}(t) =\frac{1}{\sqrt{t_{j+1}-t_{j}}}, \quad m_{j,k}(t) = \sqrt{\frac{2}{t_{j+1}-t_{j}}} \cos \left( \frac{(k-1) \, \pi \, (t-{t_{j}}) }{t_{j+1}-t_{j}} \right), \quad k \geq 2,
\end{align*}
where $t \in [t_{j},t_{j+1}]$. We utilize either a first- or a second-order time integration method.

In our numerical examples, the dynamics converge to an invariant measure. To demonstrate the convergence behavior of our algorithm, we compare our results to the second order exact statistics or Monte Carlo simulations with sufficiently high sampling rate $M_{samp}$ where the exact solution is not available. Comparisons involve the following relative pointwise errors:
\begin{align*}
\epsilon_{\mbox{mean}}(t) :=\bigg  \vert \frac{\mu_{\mbox{pce}}(t)- \mu_{\mbox{exact}}(t)}{\mu_{\mbox{exact}}(t)} \bigg \vert, \quad  \epsilon_{\mbox{var}}(t) :=\bigg  \vert \frac{\sigma^2_{\mbox{pce}}(t)- \sigma^2_{\mbox{exact}}(t)}{\sigma^2_{\mbox{exact}}(t)} \bigg \vert,
\end{align*}
where $\mu$ and $\sigma^2$ represent the mean and the variance, respectively. In some cases, we give estimations of the first six cumulants of the invariant measure. 

\begin{exmp} \rm

As a first example we consider an OU process
\begin{align} \label{eq:UOU}
 d u(t) = b_u \,(\mu_u- u(t)) \, dt + \sigma_u \, dW(t), \quad t \in [0,T], \quad u(0) = u_0,
\end{align} 
where the damping parameter is random and uniformly distributed in $[1,3]$, i.e. $b_u \sim U(1,3)$.  This is an example of 2-dimensional non-Gaussian dynamics that may be seen as a coupled system for $(u,b_u)$ with $db_u=0$.

For the first simulation, we consider the time domain $[0,4]$. The mean-return parameter $\mu_u$ is set to be $0.2$. The initial condition is normally distributed $u_0 \sim N(1,0.04) \indep W(t)$ and $\sigma_u=4$. We use the Gauss-Hermite rule for the initial condition with the level $\lambda_{u_0}=3$, whereas for the damping parameter $b_u$, we use the Gauss-Legendre quadrature rule with a varying $\lambda_{b_u}$. For Brownian motion, we use $2$-dimensional approximation with the product Gauss-Hermite rule of the level $\lambda_{\bxi}=2$. We also take the set $T_{\balp}(u)$ as the normalized Hermite polynomials. For time integration, we utilize a second-order weak Runge-Kutta  method with $ \Delta \tau=5$E-$4$.

In Figure \ref{fig:UOU}, we compare second order statistics obtained by our method to the exact solutions using $N=1,2$ and $3$, $\lambda_{b_u}=2,4,$ and $8$, and $\Delta t = 0.4,0.2,0.1,$ and $0.05$.  The results are obtained by calculating the moments in $u$ variable and then taking averaging with respect to the measure of $b_u$. The rows of the figure correspond to $N$-, $\lambda_{b_u}$-, and  $\Delta t$-convergence of the method while keeping the other two degrees of freedom constant, respectively. For each different restart step $\Delta t$, we keep the time-integration step $\Delta \tau$ the same. 

$N$-convergence of the method, Figure \ref{fig:UOUa} and Figure \ref{fig:UOUb}, shows that first two moments can be captured accurately with $N=2$. The level of accuracies are $O(10^{-9})$ and $O(10^{-5})$ for the mean and the variance, respectively. We notice that using a larger number quadrature level $\lambda_{b_u}$ and more frequent restarts also help to reduce the relative errors. We also observe that the convergence behavior of the method in terms of the size of time interval is at least quadratic; i.e. $O(\Delta t^2)$.

\begin{figure}[!htb]
\centering

\subfloat[ $\lambda_{b_u}=8$ and $\Delta  t =0.05$ ]{ \label{fig:UOUa}
\includegraphics[scale=0.4]{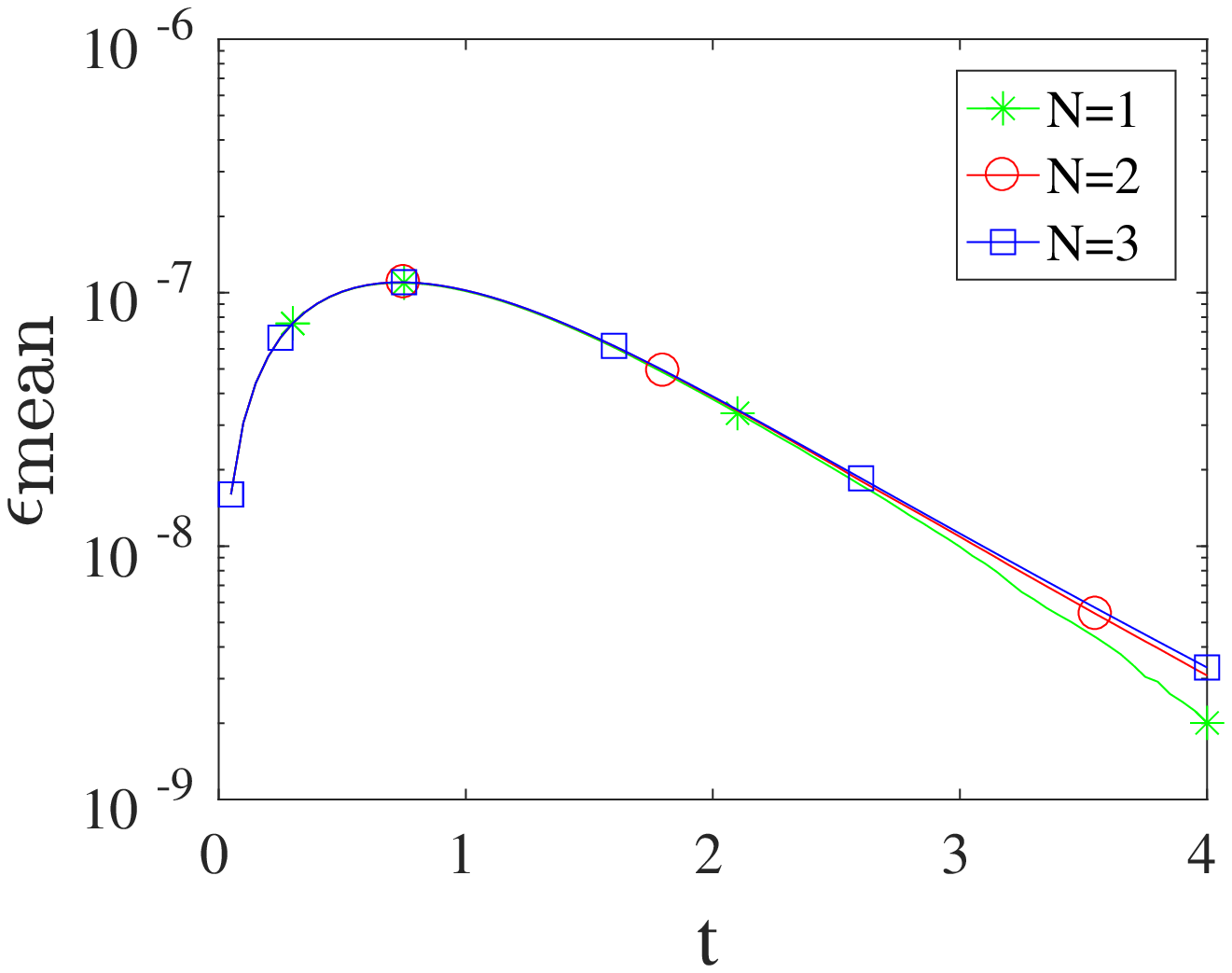}} 
\subfloat[ $\lambda_{b_u}=8$ and $\Delta  t =0.05$ ]{
\label{fig:UOUb}
\includegraphics[scale=0.4]{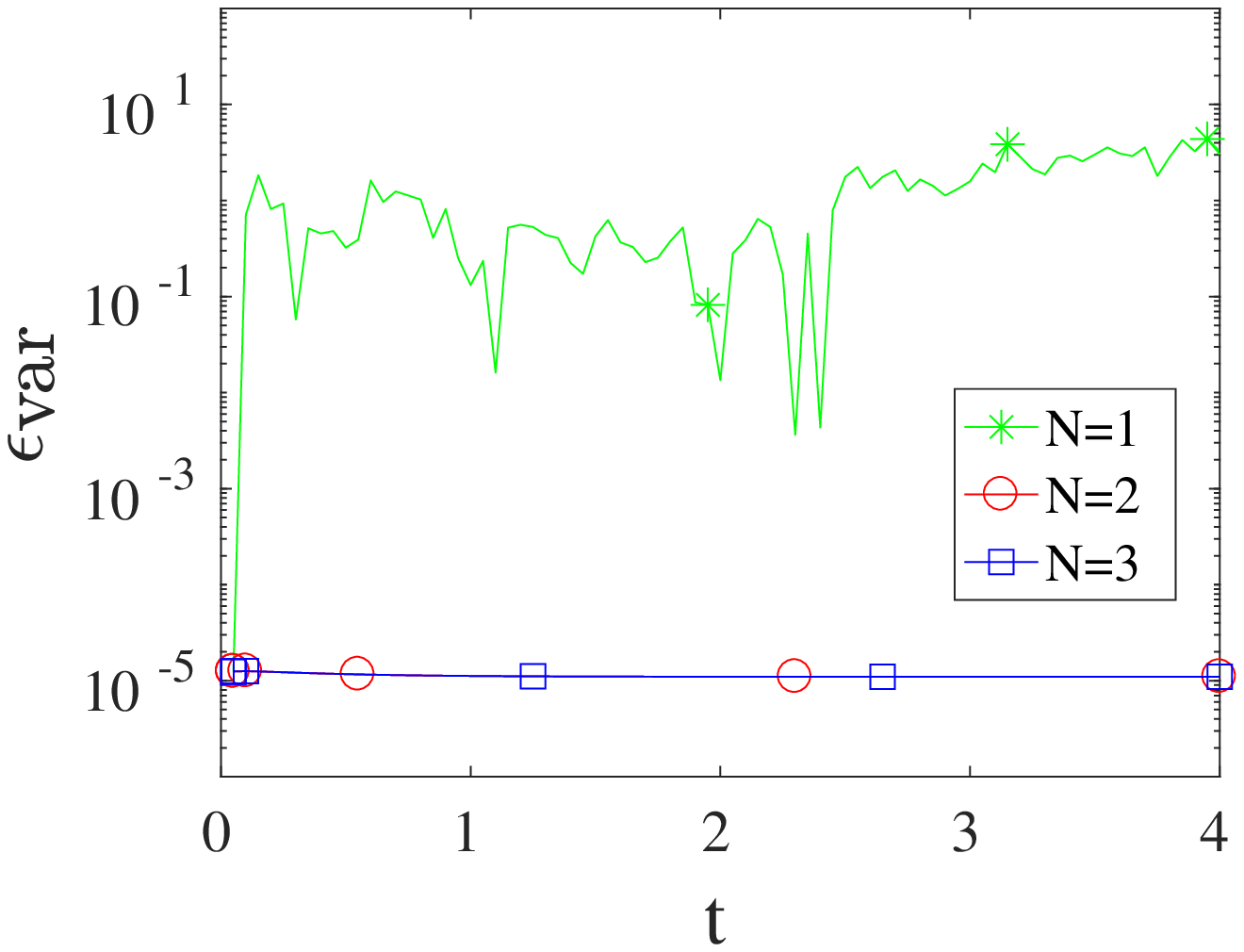}}  \\
\subfloat[ $N=2$ and $\Delta  t =0.05$ ]{
\includegraphics[scale=0.4]{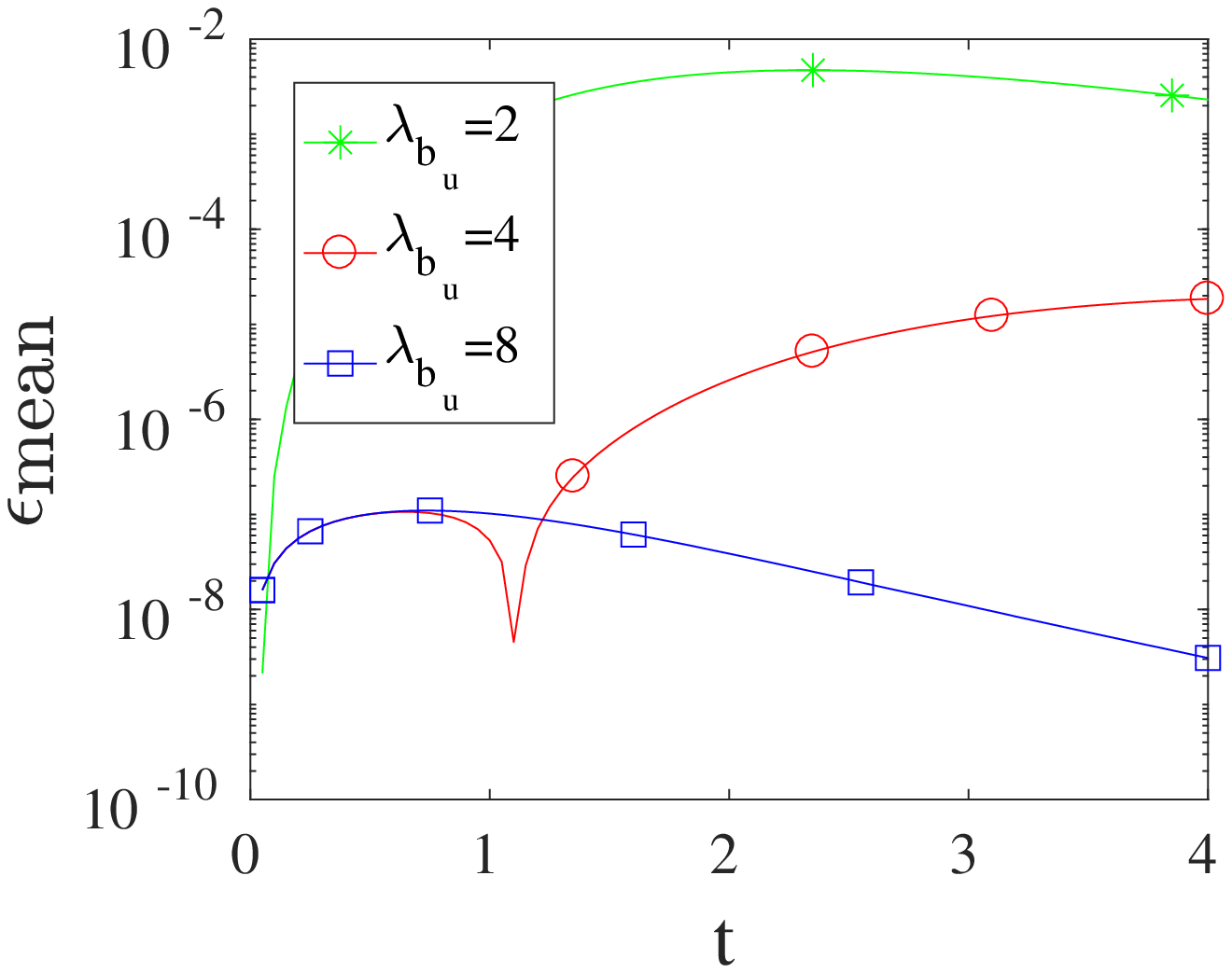}} 
\subfloat[ $N=2$ and $\Delta  t =0.05$ ]{
\includegraphics[scale=0.4]{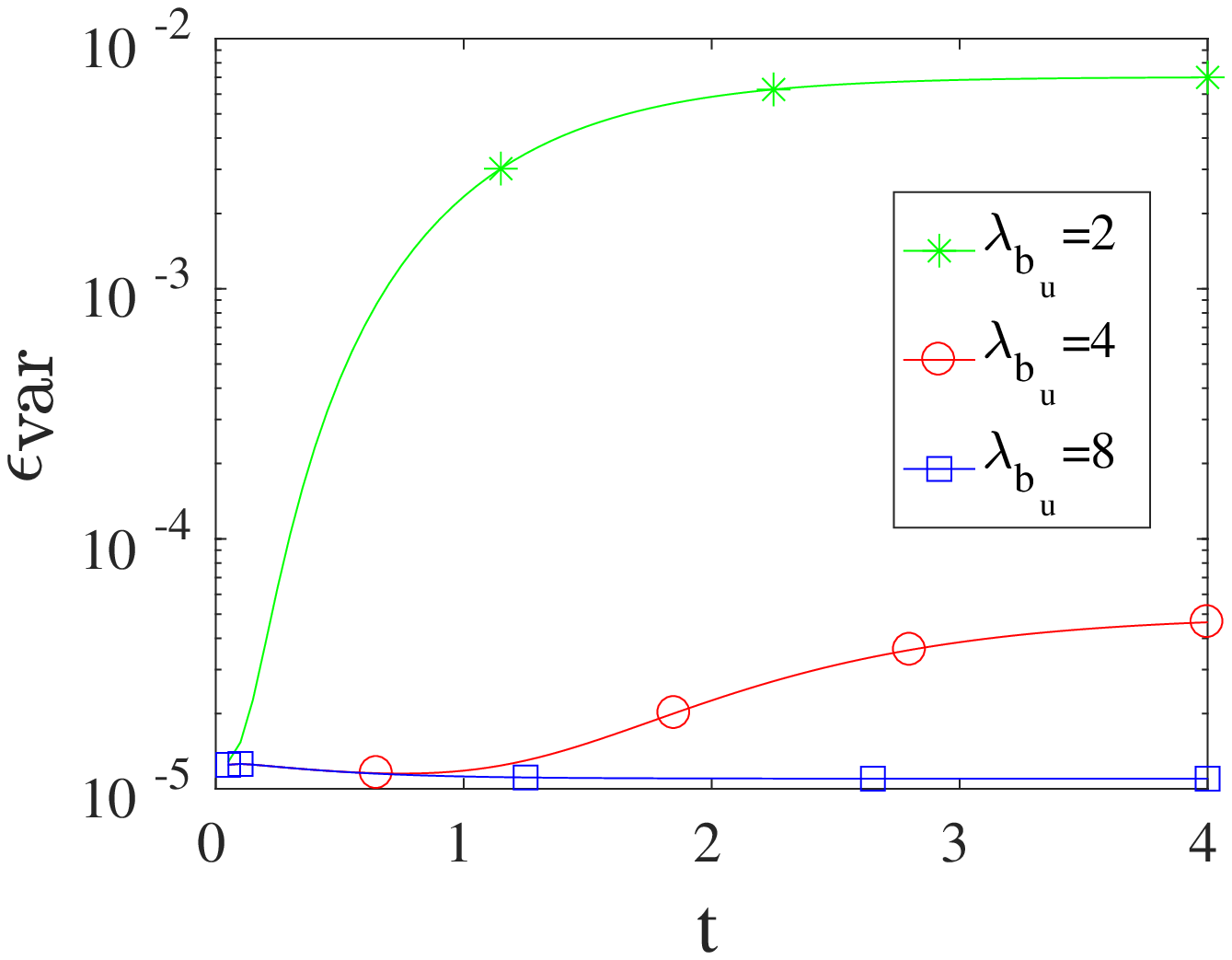}} \\
\subfloat[  $\lambda_{b_u}=8$ and $N=2$]{
\includegraphics[scale=0.4]{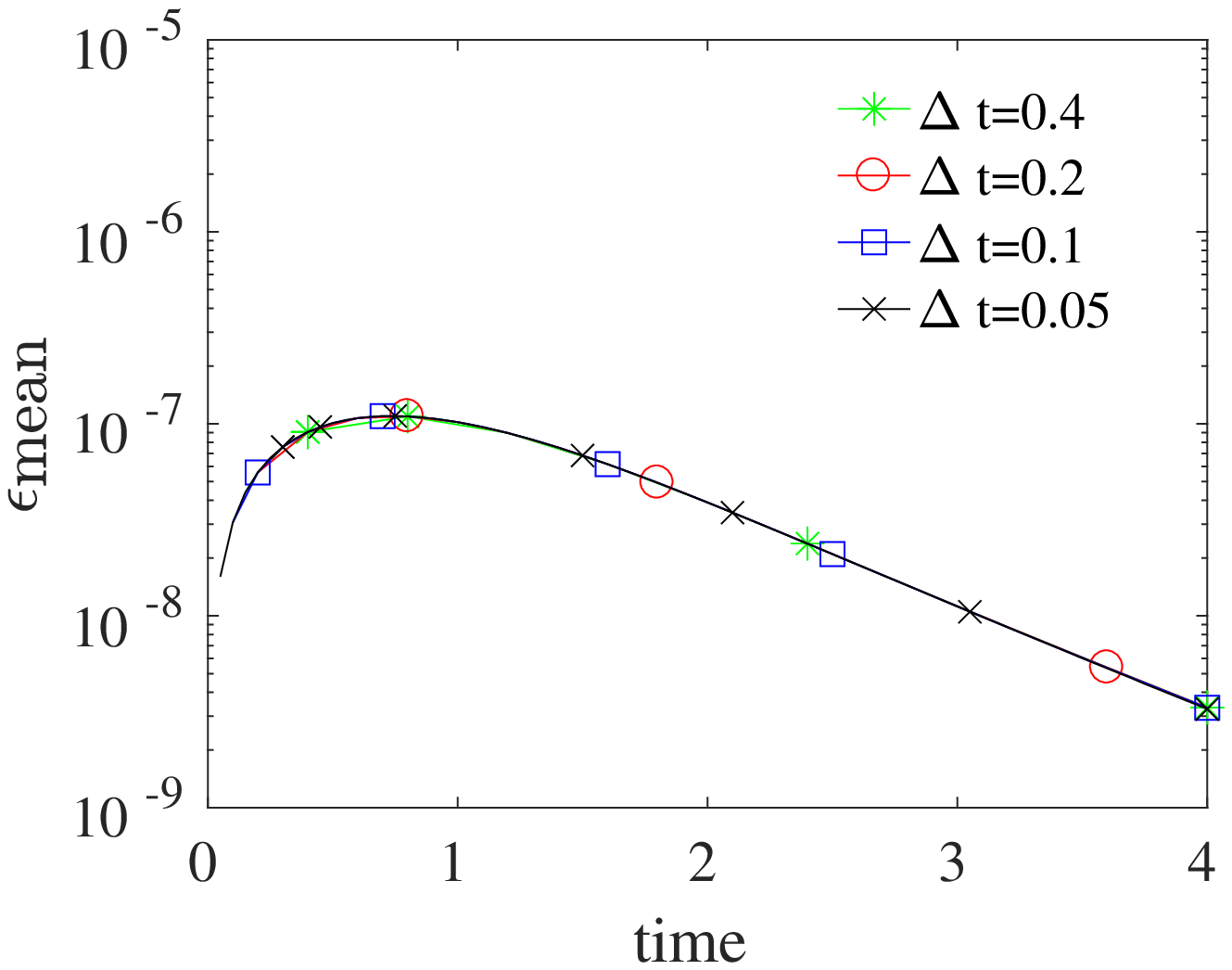}} 
\subfloat[  $\lambda_{b_u}=8$ and $N=2$ ]{
\includegraphics[scale=0.4]{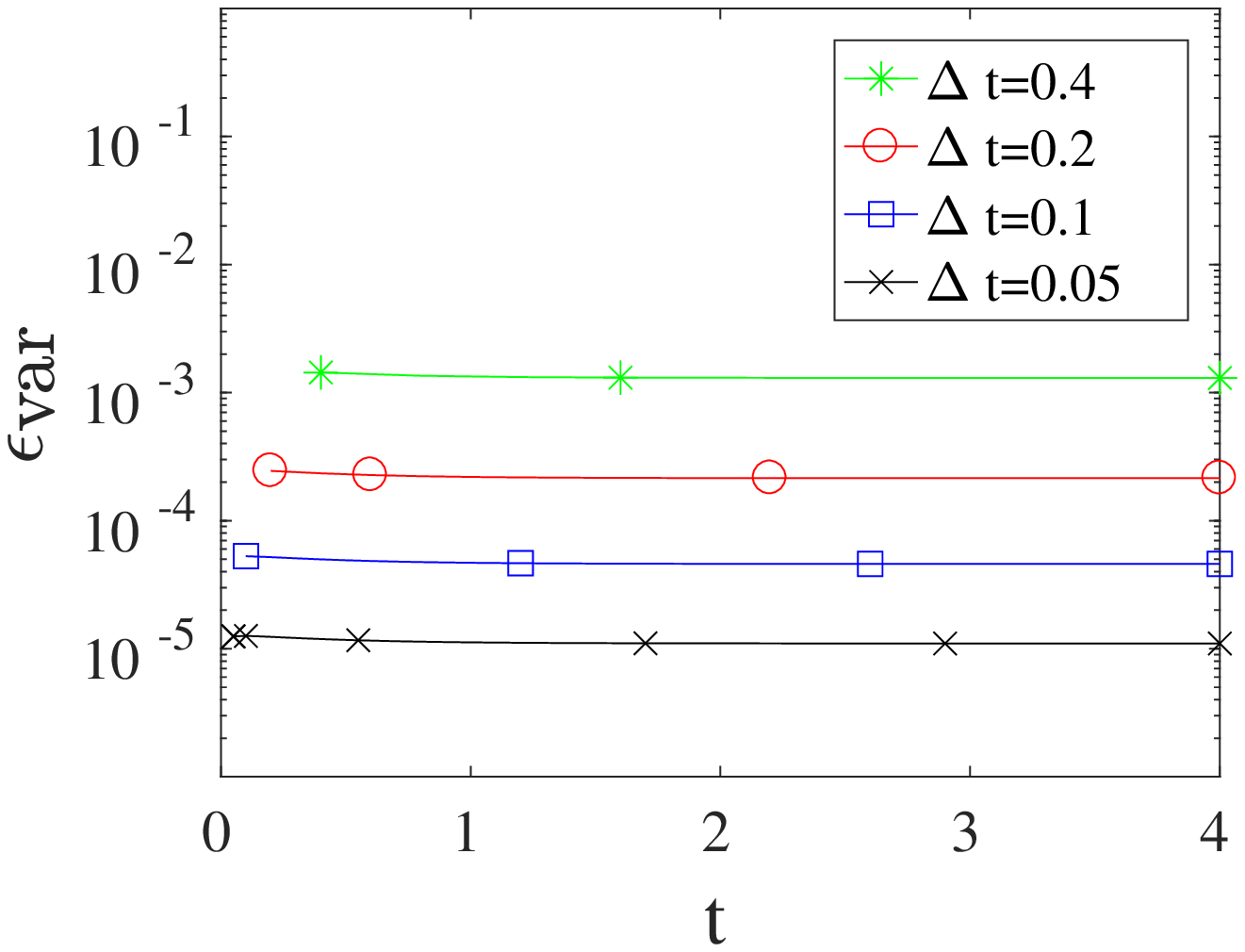}} 
\caption{Convergence behaviors in $N$, $\lambda_{b_u}$, and $\Delta t$ for the Ornstein-Uhlenbeck process with random damping.}
\label{fig:UOU}
\end{figure}

Table \ref{table:cumulants_UOU} exhibits the first six cumulants of the equation \eqref{eq:UOU} with $\mu_u=0.0$ in the long-time $T=8$. The limiting stationary measure can be obtained by solving the corresponding standard Fokker-Planck equation; see \cite{BM13, OB16, Kloeden, OB03}. The first $6$ cumulants, $\kappa_i$, $i=1,\hdots 6$, are obtained by computing moments in  the solution variable and then averaging with respect to the damping parameter. The table shows that increasing the degree $N$ of polynomials in the constraints in \eqref{eq:optL1} clearly helps to accurately capture the higher cumulants in the long-time. We also note that the approximations for the cumulants $\kappa_5$ and $\kappa_6$ become accurate when $N=6$ is used while they are inaccurate for $N=4$. This type of convergence behavior is related to the fact that the equation is linear in this case.  

\begin{table}[!htb]
\begin{minipage}[b]{1\linewidth}
\centering
\renewcommand{\tabcolsep}{0.1cm}
\renewcommand{\arraystretch}{1.00}
 \begin{tabular} { |r | r | r |r|r|r|r|}
  \hline
   & $\kappa_1$&$\kappa_2$ &$\kappa_3$  & $\kappa_4$&$\kappa_5$ & $\kappa_6$\\ \hline
   DSGC: $N=4$ & 2.09E-5 &   4.39 & 1.75E-4 & 6.06 &  -7.38 & 76.00\\ 
      DSGC: $N=6$ & 2.09E-5 &   4.39 & 1.75E-4 & 6.06 &  1.96E-3 &33.85\\ 
    Fokker-Planck &  0&     4.39 &   0 &6.06 & 0  &33.85\\ 
    \hline  
  \end{tabular}
  \caption{Cumulants obtained by Algorithm \ref{alg:DSGC} and Fokker-Planck equation at $T=8$ for the Ornstein-Uhlenbeck process with random damping.}
  \label{table:cumulants_UOU}
\end{minipage}
\end{table}

\end{exmp}

\begin{exmp} \rm 
We now consider a nonlinearity in the equation so that the damping term includes a cubic term:
\begin{align*}
 d u = - (u^2 + 1)u \,dt + \sigma_u \, dW, \quad u(0)=1.
\end{align*}
We take $T=4$, $\Delta t= 0.04$, and $\sigma_u=2$.  For Brownian motion we use $K=2$ dimensional $\bxi$ vector with a sparse Gauss-Hermite rule with the level $\lambda_{\bxi}=3$. We also take the set $T_{\balp}(u)$ as Hermite polynomials. 

In this case, we observe from Table \ref{table:cumulants_cubic} that the Algorithm \ref{alg:DSGC} can be executed accurately in the long-time by increasing the number of degrees of freedom $N$. Comparing Table \ref{table:cumulants_cubic} and Table \ref{table:cumulants_UOU}, we see that the accuracy in the higher cumulants are slightly decreased and it is harder to capture higher moments as this is a more complicated dynamics having a nonlinearity.

\begin{table}[!htb]
\begin{minipage}[b]{1\linewidth}
\centering
\renewcommand{\tabcolsep}{0.1cm}
\renewcommand{\arraystretch}{1.00}
 \begin{tabular} { |r | r | r |r|r|r|r|}
  \hline
   & $\kappa_1$&$\kappa_2$ &$\kappa_3$  & $\kappa_4$&$\kappa_5$ & $\kappa_6$\\ \hline
    DSGC: $N=4$ & 1.26E-2 &   7.35E-1 & 4.12E-2 & -3.22E-1 & 2.08E-1 &  4.75E-1\\ 
    DSGC: $N=6$ & -2.01E-3 &   7.34E-1 & 8.20E-3 & -3.39E-1 &  -7.27E-2  & 9.35E-1 \\ 
     DSGC: $N=8$ & 3.48E-4 &   7.33E-1 & -2.91E-3 & -3.39E-1 &  2.40E-3 & 9.52E-1\\ 
    Fokker-Planck &  0&     7.33E-1 &   0 & -3.39E-1 & 0  & 9.64E-1\\ 
    \hline  
  \end{tabular}
  \caption{Cumulants obtained by Algorithm \ref{alg:DSGC} and Fokker-Planck equation at $T=4$ for the cubic nonlinearity.}
  \label{table:cumulants_cubic}
\end{minipage}
\end{table}
\end{exmp}

\begin{exmp} \label{ex:CIR} \rm 
In this example, we consider a multiplicative noise term and the Cox-IngerSoll-Ross (CIR) model 
\begin{align} \label{eq:CIR}
du(t) = b_u (\mu_u-u(t)) dt + \sigma_u \sqrt{u(t)} \, dW(t),
\end{align}
which is used to describe the evolution of interest rates~\cite{CIR85}. The same model is also used in the Heston model to model random volatility. We impose the condition $2 b_u \mu_u \geq \sigma_u^2$ so that the process stays positive. The process has a stationary distribution in the long time and the second order statistics can be computed analytically. Note also that the noise amplitude is non-Lipschitz. 

We apply the first-order Milstein method \cite{Kloeden}
\begin{align} \nonumber
u(\tau_{i+1})& = u(\tau_{i}) + b_u \left(\mu_u-u(\tau_{i})-\frac{\sigma_u^2}{4 b_u} \right) dt  \\
 \label{eq:particleSys_multiplicative}
& \qquad + \sigma_u \sqrt{u(\tau_i)} (W(\tau_{i+1})-W(\tau_{i})) + \frac{\sigma_u^2}{4} (W(\tau_{i+1})-W(\tau_{i}))^2,
\end{align}
and then approximate Brownian motion increments by their finite-dimensional approximations. Since the noise term involves a function of $u$, it is necessary to use a time-integration method which takes this dependence into account. We set $T=3, u_0=1$, $\mu_u=0.6$, $\sigma_u =0.5$, $b_u=2$, $K=2$, $\lambda_{\xi}=4$, and $\Delta t=0.1$. We compute the associated orthogonal polynomials when we solve the optimization procedure.

The first row of Figure \ref{fig:DSGC_CIR} shows the evolution of the mean and the variance of the model \eqref{eq:CIR} obtained by the analytical solution and Algorithm \ref{alg:DSGC} with $N=4$. We observe that in the long time, the statistics become stationary and are captured well by the algorithm. The second row shows that the algorithm can be performed accurately in the long-time with a relatively small noise-level $\sigma_u=0.5$.

\begin{figure}[!htb]
\centering
\subfloat[$N=4$ approximation to the mean]{ \label{fig:UOUa}
\includegraphics[scale=0.4]{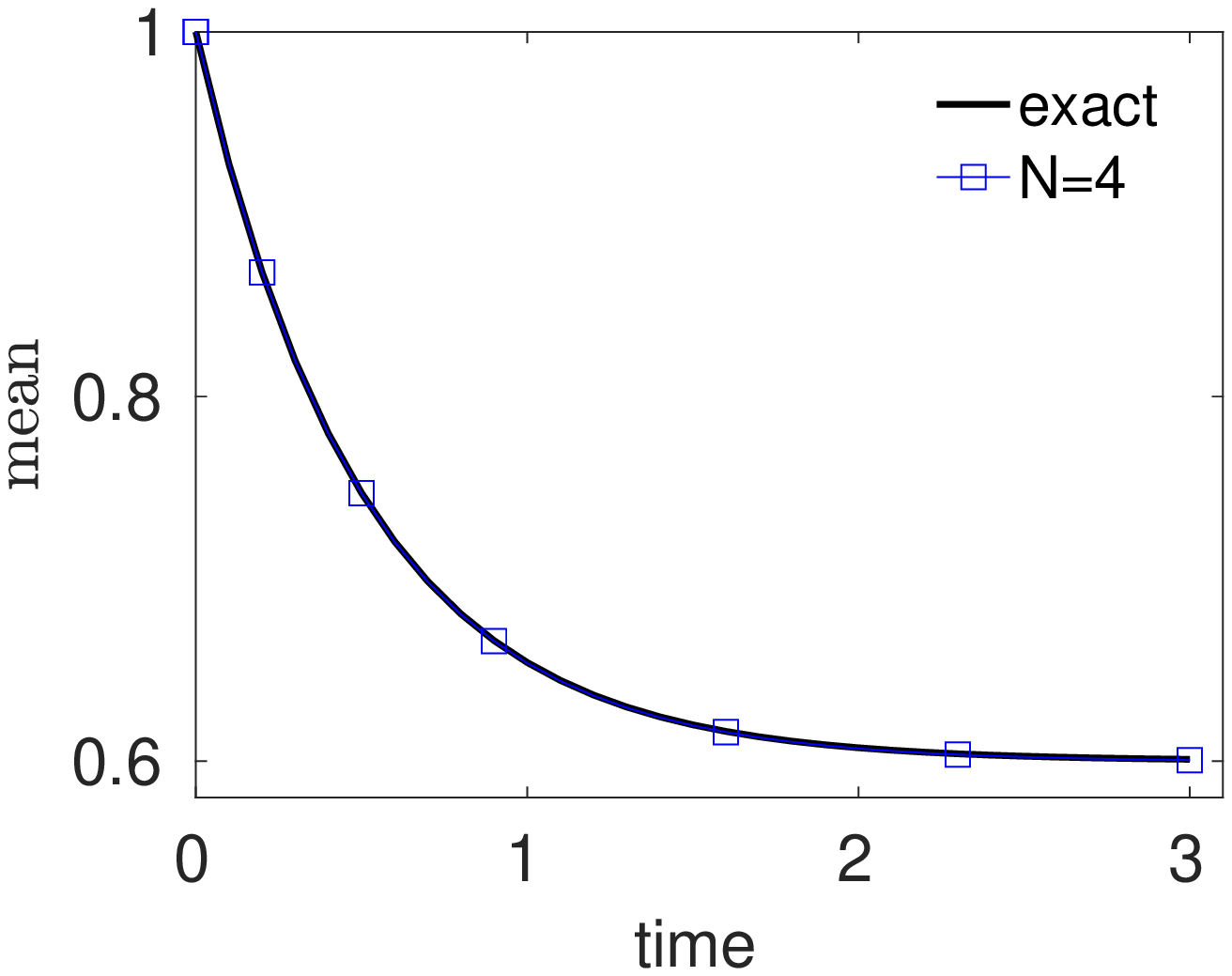}} 
\subfloat[$N=4$ approximation to the variance]{
\label{fig:UOUb}
\includegraphics[scale=0.4]{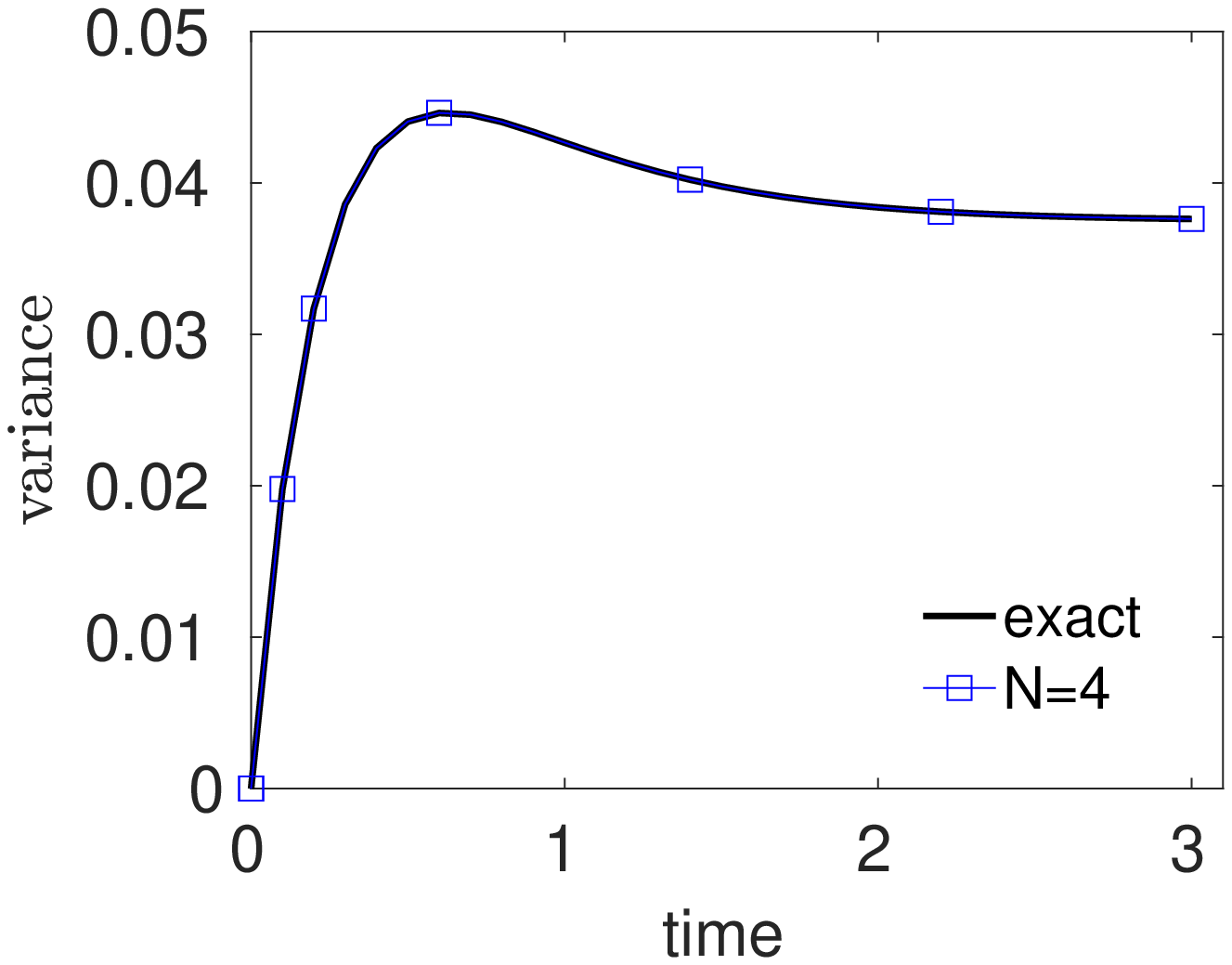}}  \\
\subfloat[ $\epsilon_{\mbox{mean}}$ for $N=1,2,4$]{
\includegraphics[scale=0.4]{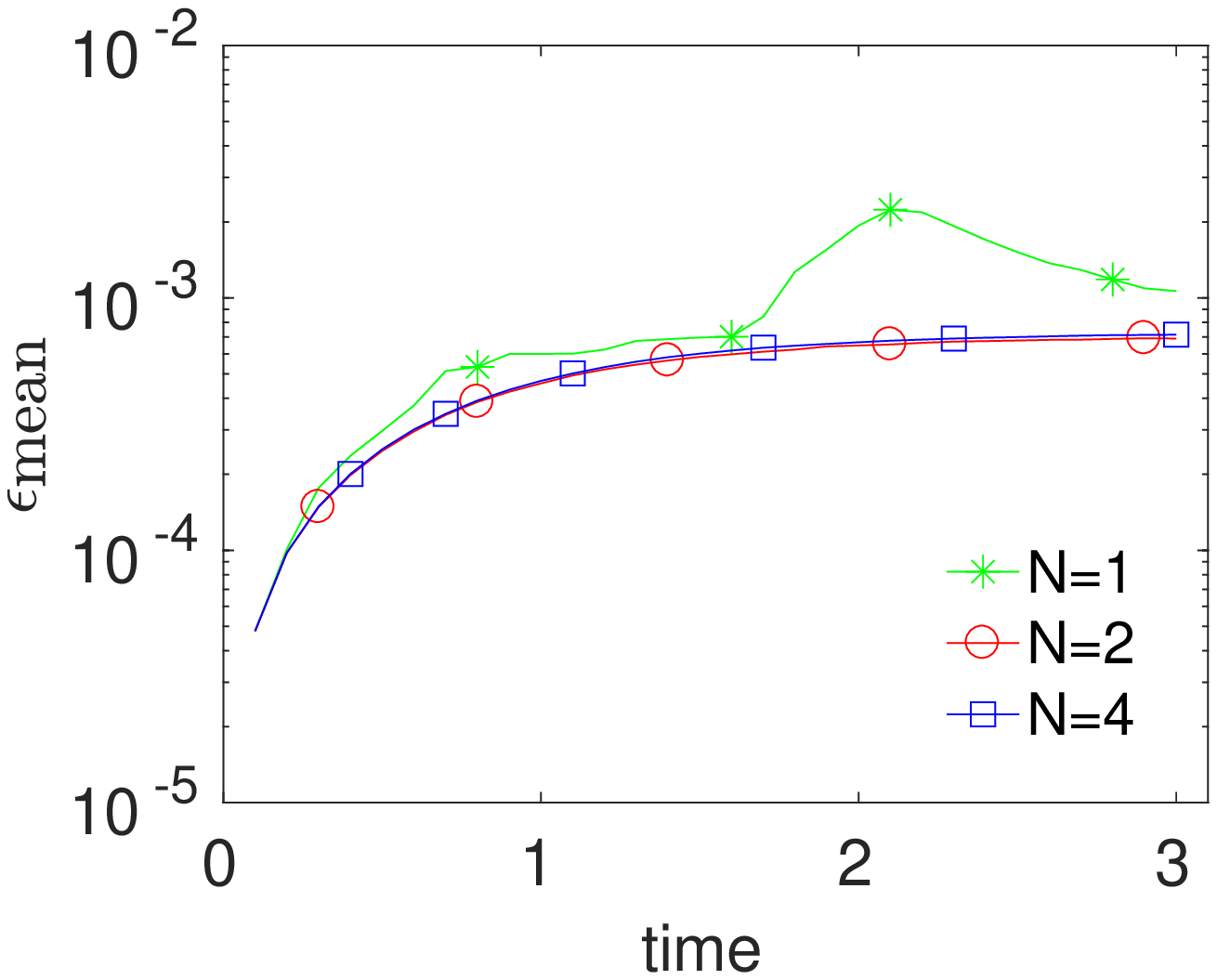}} 
\subfloat[ $\epsilon_{\mbox{var}}$ for $N=1,2,4$ ]{
\includegraphics[scale=0.4]{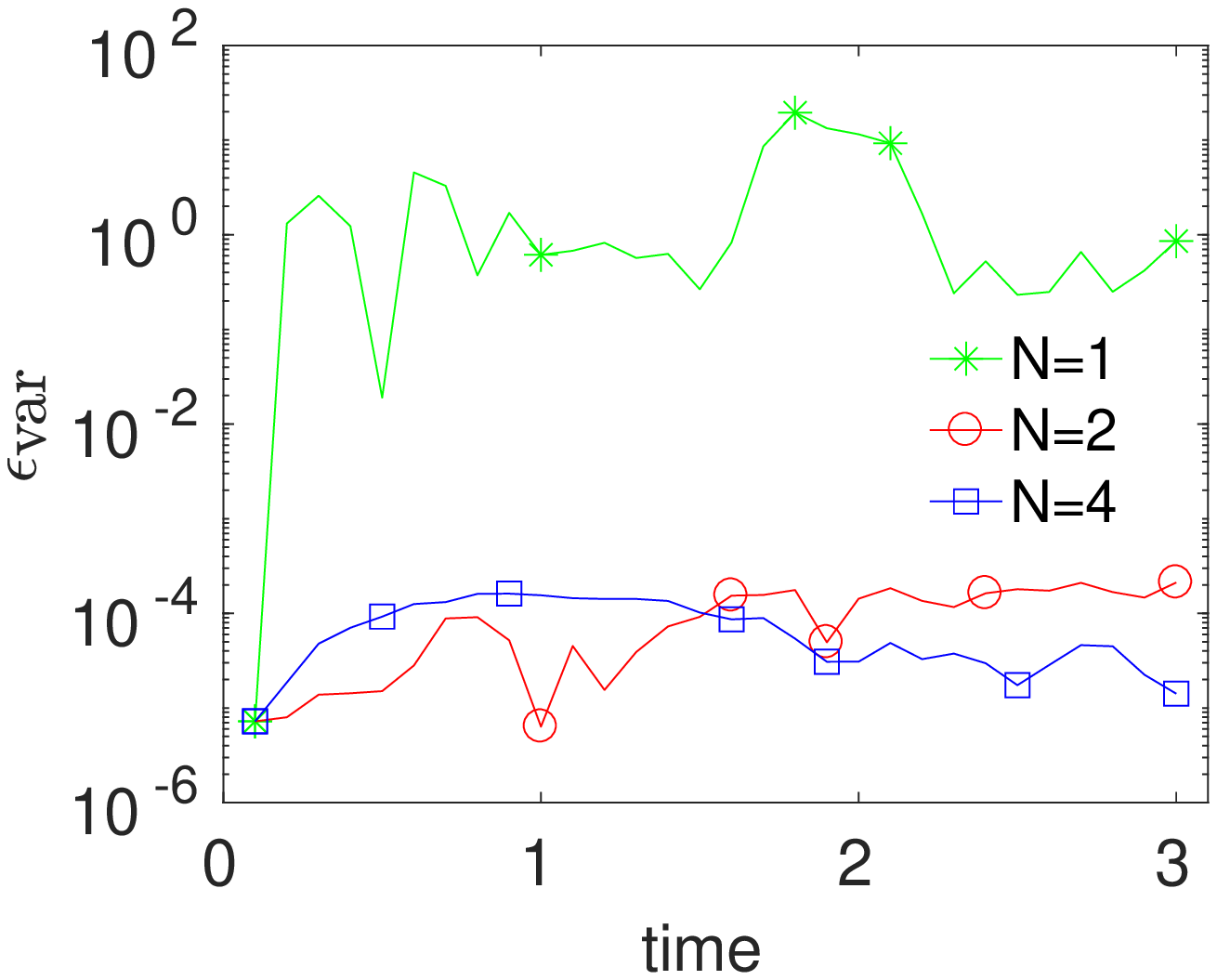}}
\caption{Evolution of the mean and the variance of the CIR model, and the convergence behaviors in $N$.}
\label{fig:DSGC_CIR}
\end{figure} 

In Table \ref{table:MC_comparison}, we compare the DSGC algorithm and MC method in terms of accuracy and timing. For a fair comparison and to minimize the error of time integration, both methods use Milstein's method with the time step $\Delta \tau = 1$E-$4$. We compare the convergence behaviors of the variance using $N$ and $M_{samp}$ for DSGC and MC, respectively. To assess the robustness of the methods, we take $b_u=4$ and consider different magnitudes of the noise: $\sigma_u=0.5,1,$ and $2$. Since MC estimates are noisy, we repeat each MC simulation $20$ times and average the errors. We first observe that the accuracy of both methods drops when $\sigma_u$ is increased. Also, if a high accuracy is needed, DSGC algorithm can be executed faster than MC in this scenario, e.g., for $\sigma_u=1$ and $\sigma_u=2$, MC should use at least $M_{samp} = 128$E+$4$ to get to the same level of accuracy of DSGC with $N=4$, which  implies that DSGC is $20$ times faster. We also see that DSGC is very efficient for small magnitude of noise $\sigma_{u}=0.5$. The elapsed times for DSGC seem to scale quadratically with $N$. Note also that the convergence of MC method is guaranteed, however, although we observe numerical convergence in DSGC, we do not have a convergence rigorous result for the approximate solution of \eqref{eq:particleSys_multiplicative}. Similar settings in \cite{PP05,PS11,LP06} establish convergence results for SDEs with sufficiently smooth coefficients, which do not apply in this case.

\begin{table}[!htb]
\centering
\renewcommand{\tabcolsep}{0.1cm}
\renewcommand{\arraystretch}{1.1}
 \begin{tabular} { |r|c | c |c| c|}
  \hline
    Algorithm       & $\sigma_u=0.5$ &$\sigma_u=1$ &$\sigma_u=2$ & Time ratio \\ \hline
  DSGC $N=3$ & 1.2E-4 & 2.2E-3 & 5.7E-3 & 0.5\\          
  DSGC $N=4$ & 2.9E-5 & 1.5E-3 & 1.9E-3 & 0.8\\ 
  DSGC $N=5$ & 1.5E-5 & 1.2E-3 & 1.3E-3 & 1.3\\ 
  DSGC $N=6$ & 1.1E-5 & 6.3E-4 & 6.4E-4 & 2.0\\ \hline
  MC   $M_{samp}=\,$1E+4   &1.21E-2 & 2.01E-2 & 2.46E-2 & 0.125\\
  MC   $M_{samp}= \,$2E+4   & 7.7E-3 & 1.30E-2 & 1.74E-2 & 0.25\\
  MC   $M_{samp}= \,$4E+4    &5.7E-3 & 1.00E-2 & 1.20E-2& 0.5\\
  MC   $M_{samp}= \,$8E+4    & 3.5E-3 & 6.03E-3 & 7.8E-3& 1.0\\ \hline
  \end{tabular}
  \caption{Errors of the variance at $T=1$, and relative timings of DSGC and MC methods using different degrees of freedom.}
  \label{table:MC_comparison}
\end{table}

\end{exmp}

\begin{exmp} \rm 
This example concerns the following $2$-dimensional nonlinear system 
\begin{equation}
\begin{aligned} \label{eq:system}
du(t) & = -(b_u+ a_u v(t)) u(t) \, dt + \sigma_u \, dW_u(t) , \\
dv(t) &= - b_v \, v(t) \, dt + \sigma_v \, dW_v(t),
\end{aligned}
\end{equation}
where $a_u \geq 0$, $b_u,b_v>0$ are damping parameters, $\sigma_u,\sigma_v>0$ are constants, and $W_u$ and $W_v$ are two real independent Brownian motions. The second equation is an Ornstein-Uhlenbeck process and it acts as a multiplicative noise in the first equation. Depending on the regime of the parameters, the dynamics of the solution $u$
exhibit intermittent non-Gaussian behavior; see detailed discussions in \cite{GHM10,BM13}.  

Following \cite{BM13}, we consider the system parameters as $
a_u=1$ , $b_u =1.2$, $b_v=0.5$,$\sigma_u =0.5$,  $\sigma_v =0.5$ and
take the initials $
u_0 = N(1,\sigma_u^2/8b_u) \, \indep \, v_0 = N(0,{\sigma_v^2/8b_v})$. In this regime, the dynamics of $u$ are characterized by unstable bursts of large amplitude~\cite{BM13}.

We consider a long-time $T=8$ and use $K=4$ variables for the random forcing terms with the quadrature level $\lambda_{\bxi}=2$. A sparse Gauss-Hermite quadrature rule is used for the random initial conditions. Since we use a small quadrature level $\lambda_{\bxi}$, we restart frequently and take $\Delta t=0.02$ in the following simulations. The ability to use small degrees of freedom with frequent restarts is one of the main advantages of the method.  

The first simulation concerns the choice of the polynomials $T_{\balp}(u,v)$. Figure \ref{fig:cond_A} shows the condition numbers of the constraint matrix $A$ in \eqref{eq:optL1} for different choices of polynomials and varying $N$. We observe that in each case, the lowest condition numbers correspond to ones which are obtained by orthogonal polynomials with respect to the joint distribution of $(u,v)$. Although the computation of orthogonal polynomials for large degree of polynomials is unstable, once the computation is carried out, it yields a well-conditioned constraint matrix for fixed $N$. Moreover, increasing the degree of polynomials $N$ leads to overall larger condition numbers, which might affect the stability of the numerical minimization procedures for large $N$.

\begin{figure}[!htb]
\centering
\subfloat[ $N=3$ ]{ 
\includegraphics[scale=0.4]{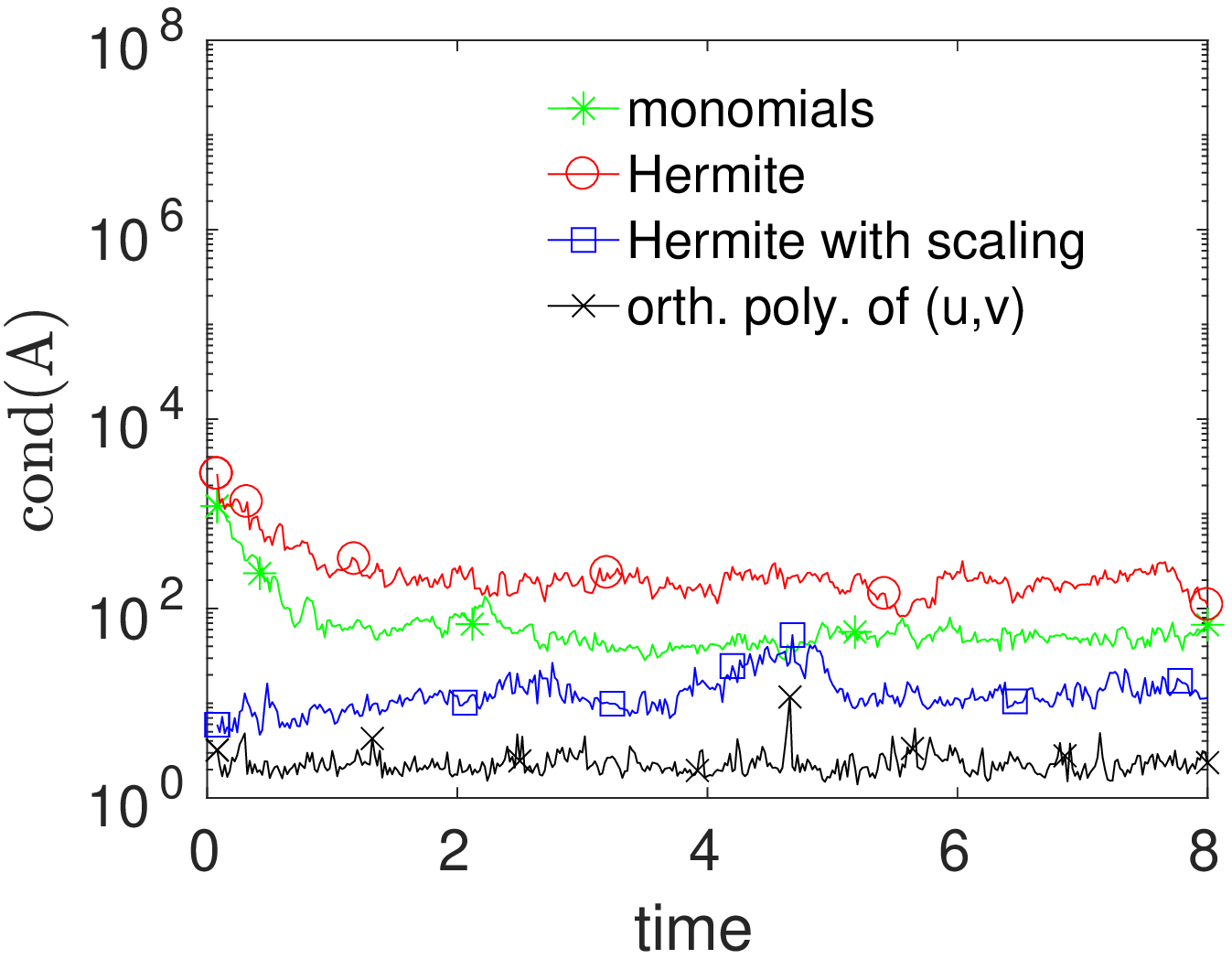}} 
\subfloat[ $N=5$ ]{
\includegraphics[scale=0.4]{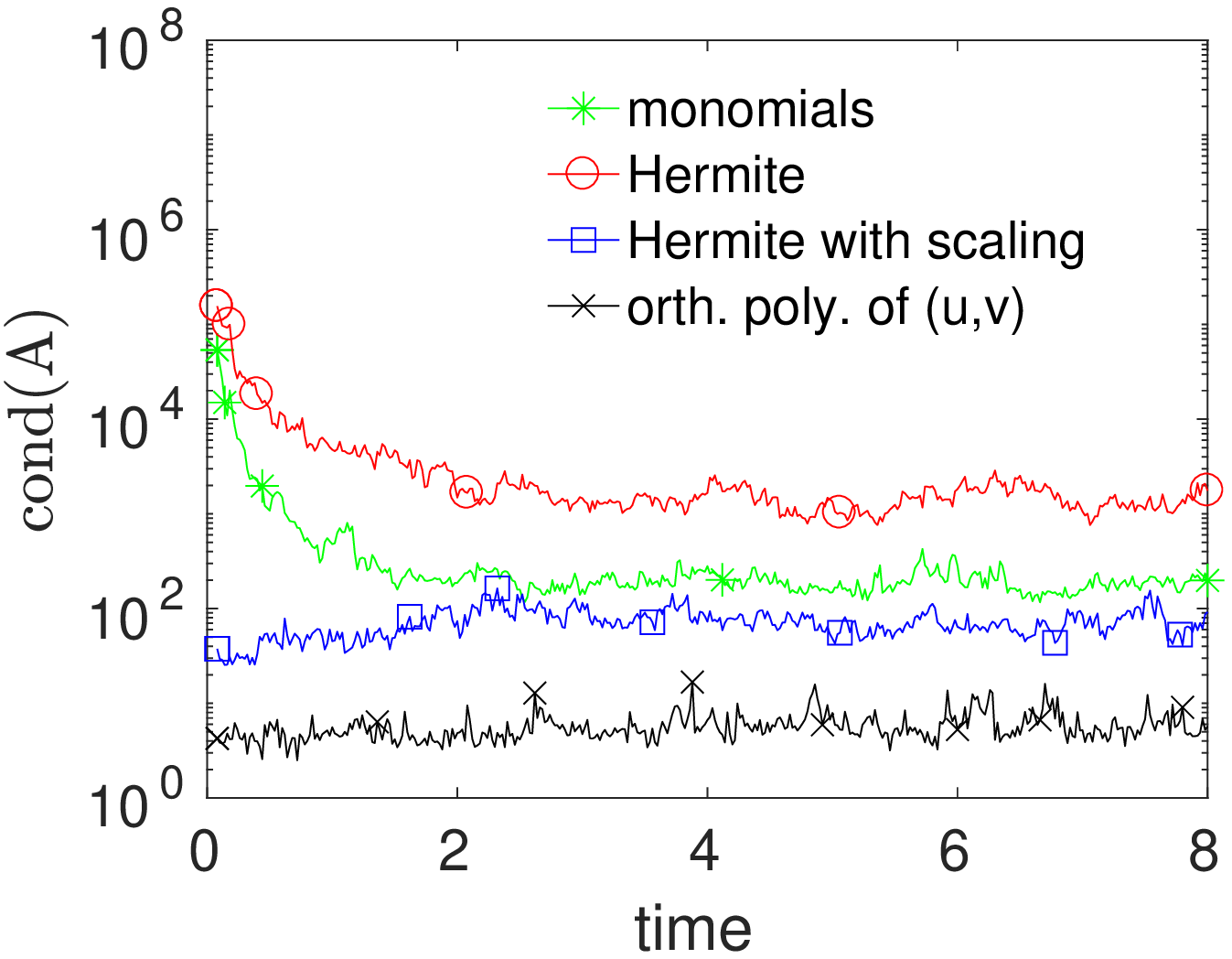}} \\ 
\subfloat[ $N=7$ ]{
\includegraphics[scale=0.4]{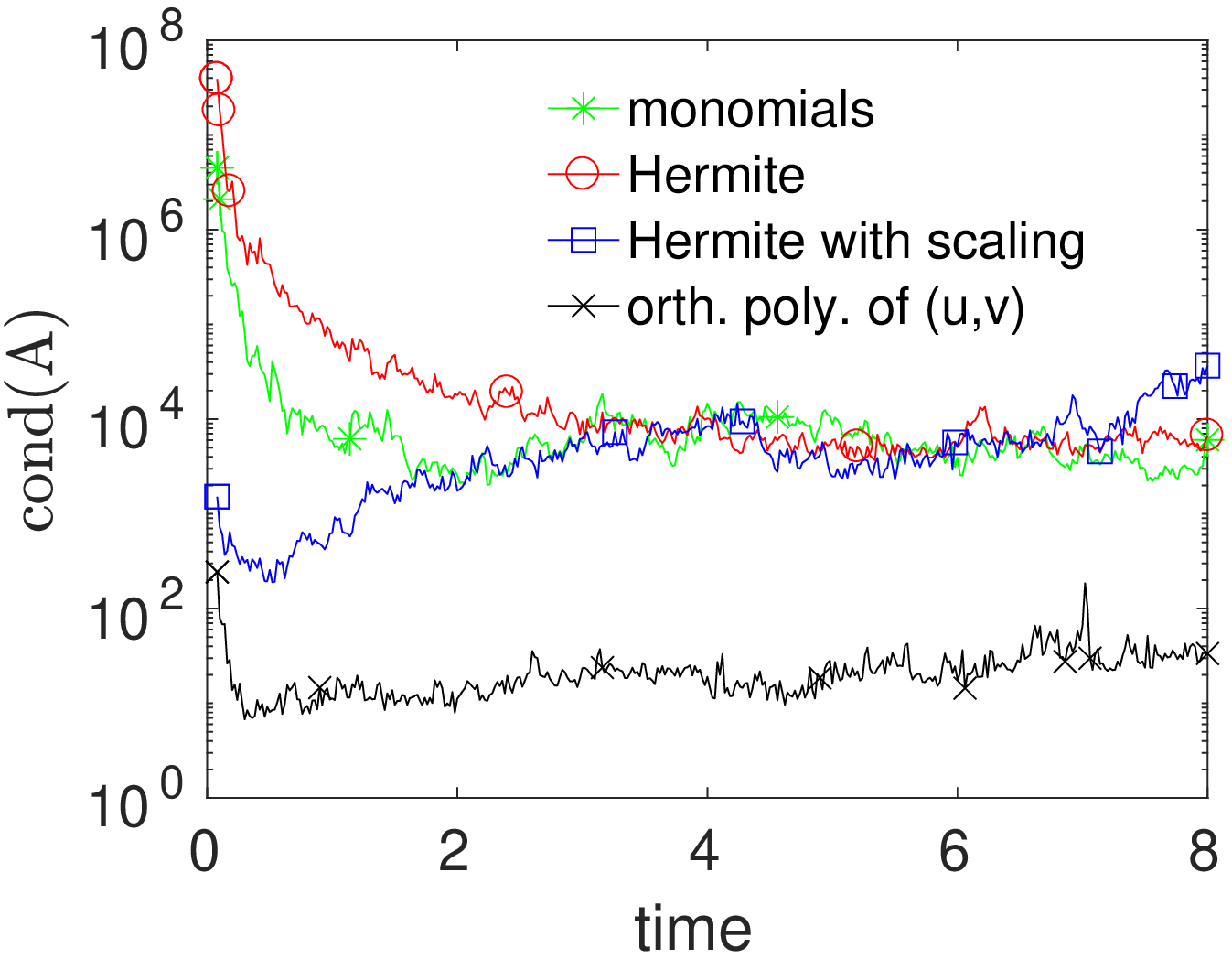}}
\caption{Condition numbers of the constraint matrix $A$ in the optimization procedure \eqref{eq:optL1} for different degrees of freedom $N$ and different choices of $T_{\balp}$. }
\label{fig:cond_A}
\end{figure} 

Using polynomials $T_{\balp}(u,v)$ which are orthogonal with respect to the joint distribution of the solution, we next demonstrate $N$-convergence of the method in terms of the mean and the variance in Figure \ref{fig:Nconv_31}. We observe from the figure that the mean and the variance can be captured in the long-time up to $O(10^{-3})$ of accuracy. Throughout the time evolution, the number of quadrature points of the joint distribution of $(u,v)$ is ${N+2 \choose 2}$. For instance, for $N=5,6,7$, and $8$, this number becomes $21,28,36$, and $45$, respectively. It is useful to note that although increasing $N$ gives better errors, it also makes the numerical optimization procedure relatively unstable. 

\begin{figure}[!htb]
\centering
\subfloat[ mean ]{ 
\includegraphics[scale=0.4]{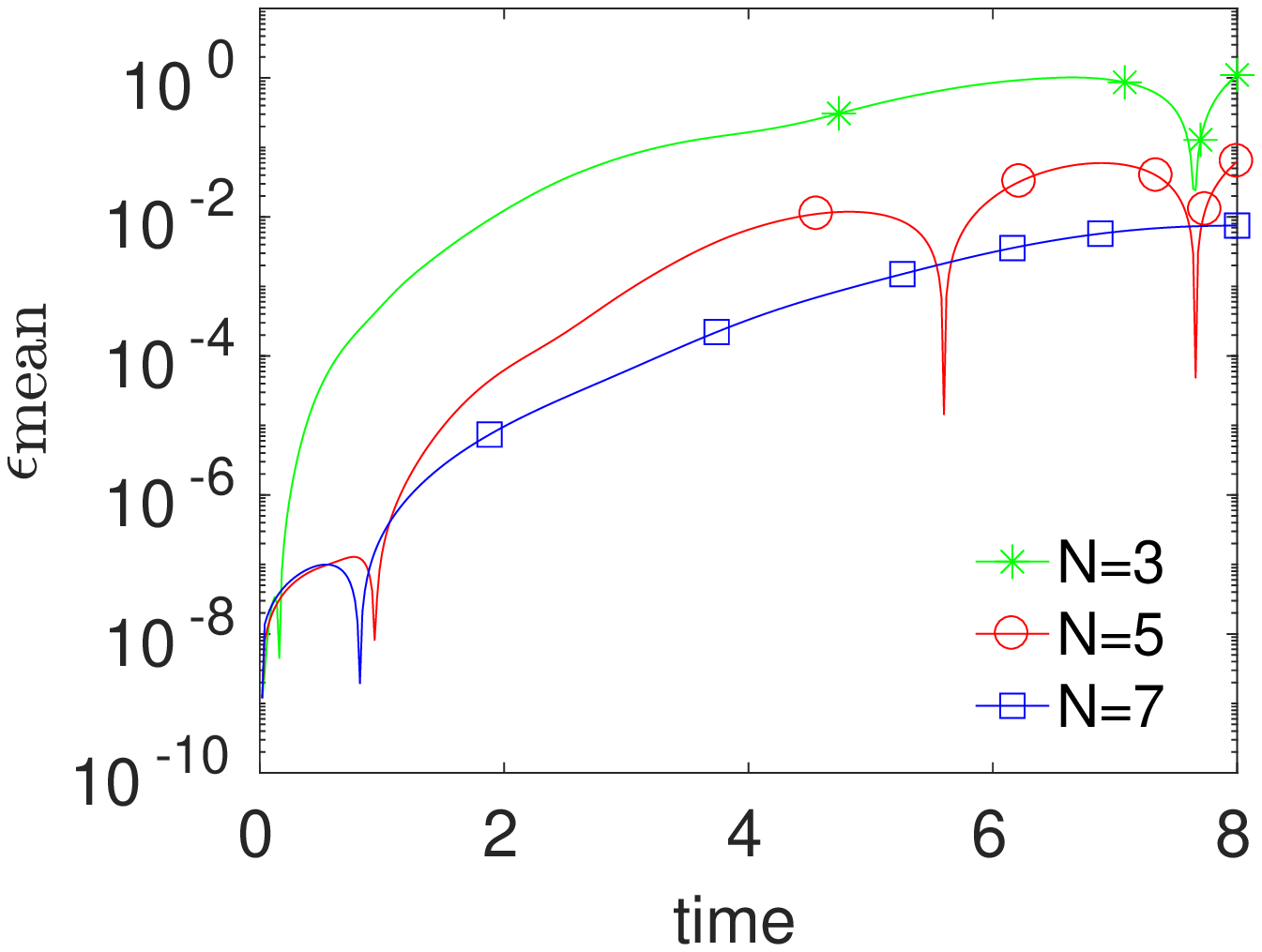}} 
\subfloat[ variance ]{ 
\includegraphics[scale=0.4]{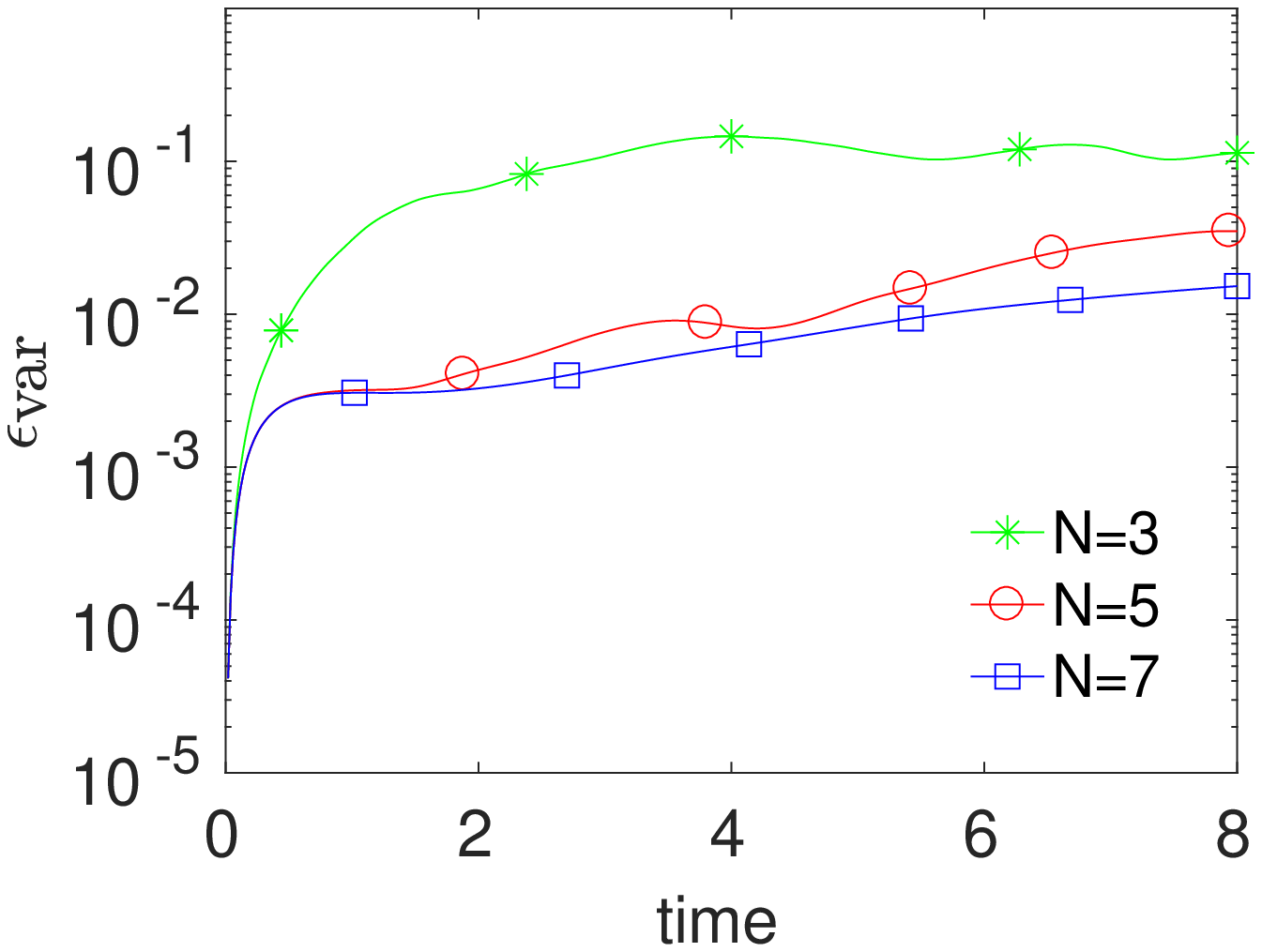}} \\ 
\subfloat[ mean ]{ 
\includegraphics[scale=0.4]{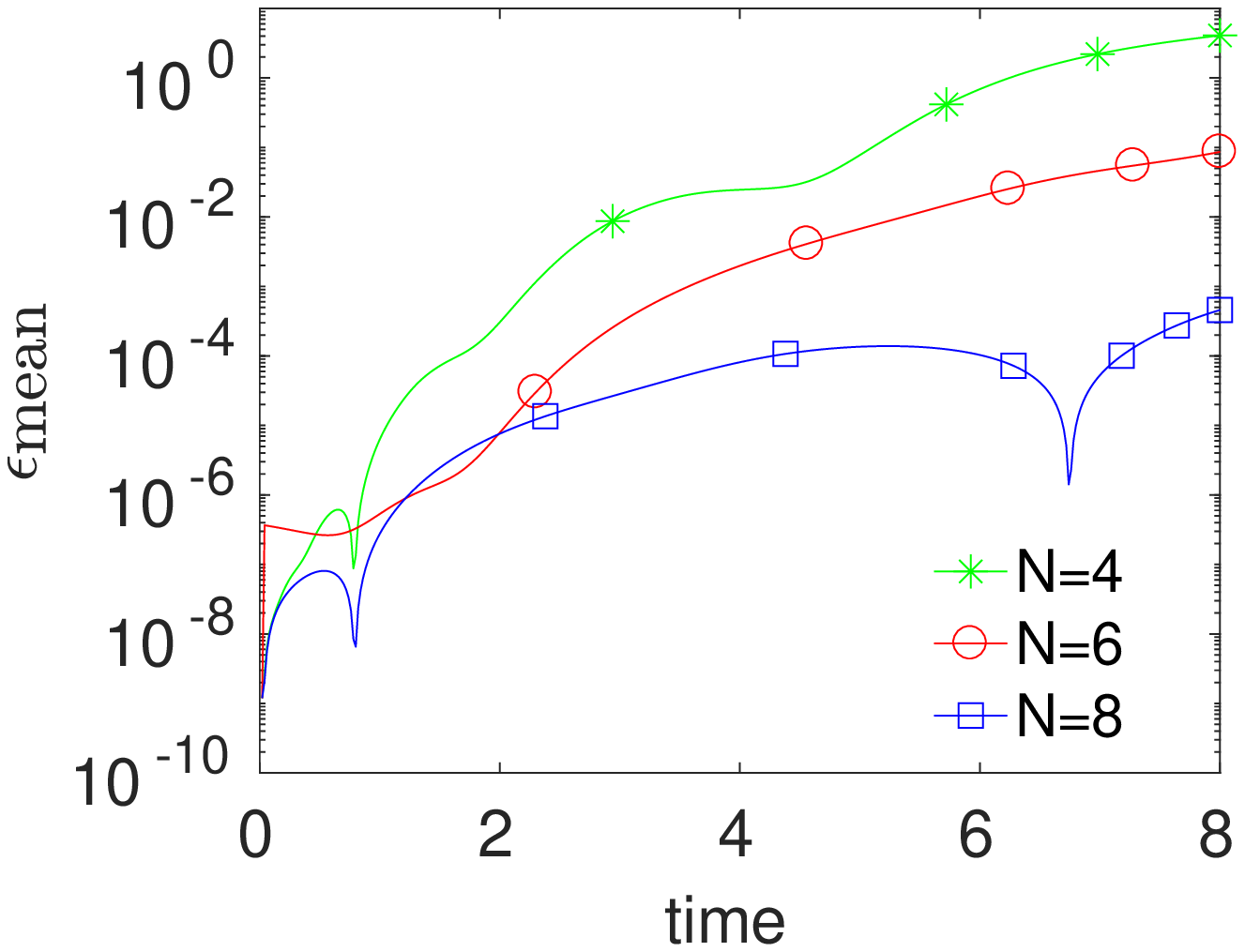}} 
\subfloat[ variance ]{ 
\includegraphics[scale=0.4]{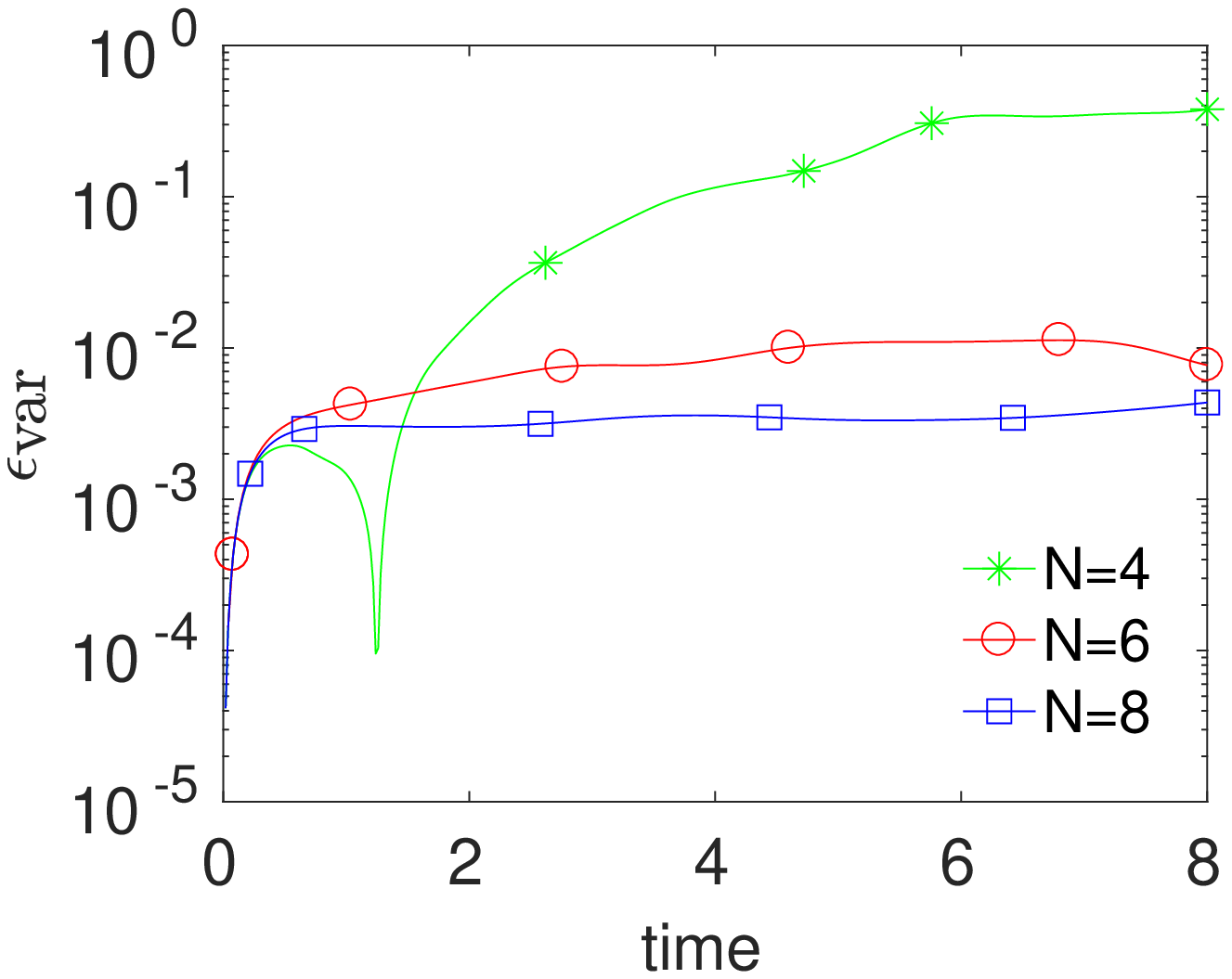}} 
\caption{$N$-convergence of the method for the nonlinear system of SDEs.}
\label{fig:Nconv_31}
\end{figure} 

Finally, we take a relatively short time $T=1$, and compare the computational times and the accuracy of DSGC, DgPC (\cite{OB17}), and MC methods. All methods use a second order time-integration method with the same time step $\Delta \tau=2$E-$4$. MC simulations are repeated 20 times to get a stable estimate of the error. Note that it is getting computationally harder to obtain a stable estimate in MC for longer times since the estimates are noisy. Both DSGC and DgPC use the restart step $\Delta t=0.05$.

From Table \ref{table:MC_DgPC_DSGC}, we first notice that for the same level of accuracy, DSGC performs better than both methods in terms of computational time. This behavior is due to the fact that the method can utilize a small number of samples and frequent restarts to achieve a high accuracy in a fast manner. We also observe that DgPC with $N=2$ offers a sufficient accuracy for many applications with a longer computational time in this scenario. Thus, DSGC offers a relatively fast way to compute moments of the solution in this low-dimensional case. Finally, we note that although DSGC is faster than DgPC in this case, it has been shown in \cite{OB17} that DgPC offers a faster way to compute the solutions of SPDEs compared to MC methods because sampling methods suffer from high number of samples and complicated dynamics. In that case, DSGC might have an advantage over MC methods since it can use small number of samples by leveraging the regularity
of the solution; see also concluding remarks about extensions to high dimensional cases.

\begin{table}[!htb]
\centering
\renewcommand{\tabcolsep}{0.1cm}
\renewcommand{\arraystretch}{1.1}
 \begin{tabular} { |r|c | c |c|}
  \hline
  Algorithm         &  $\epsilon_{\mbox{mean}}$ &  $\epsilon_{\mbox{var}}$ & Time ratio \\ \hline      
  DSGC $N=3$ & 4.6E-4 & 3.01E-2  & 0.095\\   
  DSGC $N=5$ & 9.2E-6 & 4.5E-3  &   0.2\\ 
  \hline    
  DgPC $N=2$ & 3.8E-4 & 3.7E-3  & 2.28 \\  
  DgPC $N=3$ & 3.2E-4 & 1.8E-3  & 6.71 \\ \hline  %
  MC   $M_{{samp}}= \,$4E+4    & 4.3E-3  & 6.2E-3& 0.5\\
  MC   $M_{{samp}}= \,$8E+4    & 2.6E-3  & 4.8E-3& 1.0\\ \hline
  \end{tabular}
  \caption{Errors of the mean and the variance at $T=1$, and relative timings of DSGC, DgPC, and MC methods using different degrees of freedom.}
  \label{table:MC_DgPC_DSGC}
\end{table}

\end{exmp}

%

\section{Conclusion}

We presented a collocation method, called Dynamical Sparse Grid Collocation method, which exploits the Markovian structure of SDEs to keep the degrees of freedom under control and obtain accurate long-time computations. The method uses a restart mechanism to construct a sparse quadrature rule for the solution on-the-fly and incorporates that rule with quadrature rules for the random forcing to capture the solution at later times. This fundamental idea is very similar to those of algorithms presented in our earlier works \cite{OB16, OB17} with one major difference: the current algorithm is non-intrusive. Being a non-intrusive method, one can leverage optimized legacy differential equation solvers and carry out evolution of the solution particles in parallel. 

The main computational difficulties are how to keep the number of quadrature nodes independent of time, and to compute a sparse quadrature rule for the solution variable in a stable and efficient manner on-the-fly. We made use of a $L^1$-minimization procedure with the constraints that the quadrature rule integrates polynomials in the solution variable up to certain degrees of freedom exactly. We also discussed how to extract a quadrature rule with a few nodes and presented several different polynomials bases used in the optimization. Numerical results for different nonlinear low-dimensional dynamics confirmed the ability of the algorithm reaching long times. 

For high-dimensional SDEs or SPDEs, the algorithm needs further modifications. Dynamical sparse quadrature rules for the solution will have higher number of nodes, which, in turn, may hinder the efficiency of the algorithm; especially in case of complex dynamics. Thus, further dimensionality reduction techniques and simple parallelization can be considered. Moreover, the stability of the constrained optimization in high-dimensional cases should also be investigated. We demonstrated such an extension from SDEs to SPDEs in our previous works \cite{OB16,OB17}. The DgPC method was coupled with Karhunen--Loeve expansion and a few dominating random modes at each restart were selected and incorporated into a PCE to represent future solutions. A similar extension can be done for DSGC by using the KLE of random solutions in time and constructing sparse quadrature rules for the selected finite number of KLE random modes.


%

\bibliographystyle{plain}

\begin{thebibliography}{10}

\bibitem{A10}
Rafail~V. Abramov.
\newblock The multidimensional maximum entropy moment problem: a review on
  numerical methods.
\newblock {\em Commun. Math. Sci.}, 8(2):377--392, 2010.

\bibitem{AGPRH12}
M.~Arnst, R.~Ghanem, E.~Phipps, and J.~Red-Horse.
\newblock Measure transformation and efficient quadrature in
  reduced-dimensional stochastic modeling of coupled problems.
\newblock {\em Internat. J. Numer. Methods Engrg.}, 92(12):1044--1080, 2012.

\bibitem{AGPRH14}
M.~Arnst, R.~Ghanem, E.~Phipps, and J.~Red-Horse.
\newblock Reduced chaos expansions with random coefficients in
  reduced-dimensional stochastic modeling of coupled problems.
\newblock {\em Internat. J. Numer. Methods Engrg.}, 97(5):352--376, 2014.

\bibitem{BNT07}
Ivo Babuska, Fabio Nobile, and Ra\'ul Tempone.
\newblock A stochastic collocation method for elliptic partial differential
  equations with random input data.
\newblock {\em SIAM J. Numer. Anal.}, 45(3):1005--1034, 2007.

\bibitem{BS11}
G{\'e}raud Blatman and Bruno Sudret.
\newblock Adaptive sparse polynomial chaos expansion based on least angle
  regression.
\newblock {\em J. Comput. Phys.}, 230(6):2345--2367, 2011.

\bibitem{BM13}
M.~Branicki and A.~J. Majda.
\newblock Fundamental limitations of polynomial chaos for uncertainty
  quantification in systems with intermittent instabilities.
\newblock {\em Commun. Math. Sci.}, 11(1):55--103, 2013.

\bibitem{CIR85}
John~C. Cox, Jonathan~E. Ingersoll, Jr., and Stephen~A. Ross.
\newblock A theory of the term structure of interest rates.
\newblock {\em Econometrica}, 53(2):385--407, 1985.

\bibitem{D67}
Philip~J. Davis.
\newblock A construction of nonnegative approximate quadratures.
\newblock {\em Math. Comp.}, 21:578--582, 1967.

\bibitem{DR84}
Philip~J. Davis and Philip Rabinowitz.
\newblock {\em Methods of numerical integration}.
\newblock Computer Science and Applied Mathematics. Academic Press, Inc.,
  Orlando, FL, second edition, 1984.

\bibitem{DNPKGL04}
B.~J. Debusschere, H.~N. Najm, P.~P. P{\'e}bay, O.~M. Knio, R.~G. Ghanem, and
  O.~P. Le~Ma{\^{\i}}tre.
\newblock Numerical challenges in the use of polynomial chaos representations
  for stochastic processes.
\newblock {\em SIAM J. Sci. Comput.}, 26(2):698--719, 2004.

\bibitem{gautschi}
W.~Gautschi.
\newblock {\em Orthogonal polynomials: computation and approximation}.
\newblock Numerical Mathematics and Scientific Computation. Oxford University
  Press, New York, 2004.
\newblock {O}xford Science Publications.

\bibitem{GSVK10}
M.~Gerritsma, J-B. van~der Steen, P.~Vos, and G.~E. Karniadakis.
\newblock Time-dependent generalized polynomial chaos.
\newblock {\em J. Comput. Phys.}, 229(22):8333--8363, 2010.

\bibitem{GHM10}
B.~Gershgorin, J.~Harlim, and A.~J. Majda.
\newblock Test models for improving filtering with model errors through
  stochastic parameter estimation.
\newblock {\em J. Comput. Phys.}, 229(1):1--31, 2010.

\bibitem{Ger07}
Thomas Gerstner.
\newblock {\em {S}parse {G}rid {Q}uadrature {M}ethods for {C}omputational
  {F}inance}.
\newblock Uni- versity of Bonn, 2007.

\bibitem{GG98}
Thomas Gerstner and Michael Griebel.
\newblock Numerical integration using sparse grids.
\newblock {\em Numer. Algorithms}, 18(3-4):209--232, 1998.

\bibitem{G08}
Michael~B. Giles.
\newblock Multilevel {M}onte {C}arlo path simulation.
\newblock {\em Oper. Res.}, 56(3):607--617, 2008.

\bibitem{GW69}
Gene~H. Golub and John~H. Welsch.
\newblock Calculation of {G}auss quadrature rules.
\newblock {\em Math. Comp. 23 (1969), 221-230; addendum, ibid.}, 23(106, loose
  microfiche suppl):A1--A10, 1969.

\bibitem{gb08}
Michael Grant and Stephen Boyd.
\newblock Graph implementations for nonsmooth convex programs.
\newblock In V.~Blondel, S.~Boyd, and H.~Kimura, editors, {\em Recent Advances
  in Learning and Control}, Lecture Notes in Control and Information Sciences,
  pages 95--110. Springer-Verlag Limited, 2008.
\newblock \url{http://stanford.edu/~boyd/graph_dcp.html}.

\bibitem{cvx}
Michael Grant and Stephen Boyd.
\newblock {CVX}: Matlab software for disciplined convex programming, version
  2.1.
\newblock \url{http://cvxr.com/cvx}, March 2014.

\bibitem{H01}
Desmond~J. Higham.
\newblock An algorithmic introduction to numerical simulation of stochastic
  differential equations.
\newblock {\em SIAM Rev.}, 43(3):525--546, 2001.

\bibitem{H15}
Desmond~J. Higham.
\newblock An introduction to multilevel {M}onte {C}arlo for option valuation.
\newblock {\em Int. J. Comput. Math.}, 92(12):2347--2360, 2015.

\bibitem{HLRZ06}
T.~Y. Hou, W.~Luo, B.~Rozovskii, and H-M. Zhou.
\newblock Wiener chaos expansions and numerical solutions of randomly forced
  equations of fluid mechanics.
\newblock {\em J. Comput. Phys.}, 216(2):687--706, 2006.

\bibitem{JK09}
A.~Jentzen and P.~E. Kloeden.
\newblock The numerical approximation of stochastic partial differential
  equations.
\newblock {\em Milan J. Math.}, 77:205--244, 2009.

\bibitem{Kloeden}
P.~E. Kloeden and E.~Platen.
\newblock {\em Numerical solution of stochastic differential equations},
  volume~23 of {\em Applications of Mathematics (New York)}.
\newblock Springer-Verlag, Berlin, 1992.

\bibitem{LeMK10}
O.~P. Le~Ma{\^{\i}}tre and O.~M. Knio.
\newblock {\em Spectral methods for uncertainty quantification}.
\newblock Scientific Computation. Springer, New York, 2010.

\bibitem{LL12}
C.~Litterer and T.~Lyons.
\newblock High order recombination and an application to cubature on {W}iener
  space.
\newblock {\em Ann. Appl. Probab.}, 22(4):1301--1327, 2012.

\bibitem{LMR97}
Sergey Lototsky, Remigijus Mikulevicius, and Boris~L. Rozovskii.
\newblock Nonlinear filtering revisited: a spectral approach.
\newblock {\em SIAM J. Control Optim.}, 35(2):435--461, 1997.

\bibitem{Luo}
W.~Luo.
\newblock {\em Wiener {C}haos {E}xpansion and {N}umerical {S}olutions of
  {S}tochastic {P}artial {D}ifferential {E}quations}.
\newblock ProQuest LLC, Ann Arbor, MI, 2006.
\newblock California Institute of Technology.

\bibitem{LP06}
Harald Luschgy and Gilles Pag\`es.
\newblock Functional quantization of a class of {B}rownian diffusions: a
  constructive approach.
\newblock {\em Stochastic Process. Appl.}, 116(2):310--336, 2006.

\bibitem{LV04}
Terry Lyons and Nicolas Victoir.
\newblock Cubature on {W}iener space.
\newblock {\em Proc. R. Soc. Lond. Ser. A Math. Phys. Eng. Sci.},
  460(2041):169--198, 2004.
\newblock Stochastic analysis with applications to mathematical finance.

\bibitem{MNGK04}
O.~P. Ma{\^{\i}}tre, O.~M. Knio, H.~N. Najm, and R.~G. Ghanem.
\newblock Uncertainty propagation using {W}iener-{H}aar expansions.
\newblock {\em J. Comput. Phys.}, 197(1):28--57, 2004.

\bibitem{MR04}
R.~Mikulevicius and B.~L. Rozovskii.
\newblock Stochastic {N}avier-{S}tokes equations for turbulent flows.
\newblock {\em SIAM J. Math. Anal.}, 35(5):1250--1310, 2004.

\bibitem{MT97}
G.~N. Milstein and M.~V. Tret$\prime$yakov.
\newblock Numerical methods in the weak sense for stochastic differential
  equations with small noise.
\newblock {\em SIAM J. Numer. Anal.}, 34(6):2142--2167, 1997.

\bibitem{NTW082}
F.~Nobile, R.~Tempone, and C.~G. Webster.
\newblock An anisotropic sparse grid stochastic collocation method for partial
  differential equations with random input data.
\newblock {\em SIAM J. Numer. Anal.}, 46(5):2411--2442, 2008.

\bibitem{NTW08}
F.~Nobile, R.~Tempone, and C.~G. Webster.
\newblock A sparse grid stochastic collocation method for partial differential
  equations with random input data.
\newblock {\em SIAM J. Numer. Anal.}, 46(5):2309--2345, 2008.

\bibitem{N06}
Jorge Nocedal and Stephen~J. Wright.
\newblock {\em Numerical optimization}.
\newblock Springer Series in Operations Research and Financial Engineering.
  Springer, New York, second edition, 2006.

\bibitem{NR99}
Erich Novak and Klaus Ritter.
\newblock Simple cubature formulas with high polynomial exactness.
\newblock {\em Constr. Approx.}, 15(4):499--522, 1999.

\bibitem{OB03}
B.~{\O}ksendal.
\newblock {\em Stochastic differential equations}.
\newblock Universitext. Springer-Verlag, Berlin, sixth edition, 2003.
\newblock An introduction with applications.

\bibitem{OB16}
H.~Cagan Ozen and Guillaume Bal.
\newblock Dynamical polynomial chaos expansions and long time evolution of
  differential equations with random forcing.
\newblock {\em SIAM/ASA J. Uncertain. Quantif.}, 4(1):609--635, 2016.

\bibitem{OB17}
H.~Cagan Ozen and Guillaume Bal.
\newblock A dynamical polynomial chaos approach for long-time evolution of
  {SPDEs}.
\newblock {\em Journal of Computational Physics}, 343:300 -- 323, 2017.

\bibitem{PP05}
Gilles Pag\`es and Jacques Printems.
\newblock Functional quantization for numerics with an application to option
  pricing.
\newblock {\em Monte Carlo Methods Appl.}, 11(4):407--446, 2005.

\bibitem{PS11}
Gilles Pag\`es and Afef Sellami.
\newblock Convergence of multi-dimensional quantized {SDE}'s.
\newblock In {\em S\'eminaire de {P}robabilit\'es {XLIII}}, volume 2006 of {\em
  Lecture Notes in Math.}, pages 269--307. Springer, Berlin, 2011.

\bibitem{Putinar}
Mihai Putinar.
\newblock A note on {T}chakaloff's theorem.
\newblock {\em Proc. Amer. Math. Soc.}, 125(8):2409--2414, 1997.

\bibitem{EB2015}
Ernest~K. Ryu and Stephen~P. Boyd.
\newblock Extensions of {G}auss quadrature via linear programming.
\newblock {\em Found. Comput. Math.}, 15(4):953--971, 2015.

\bibitem{smolyak}
S.~A. Smolyak.
\newblock Quadrature and interpolation formulas for tensor products of certain
  classes of functions.
\newblock {\em Soviet Math. Dokl.}, 4:240--243, 1963.

\bibitem{Tchakaloff}
Vladimir Tchakaloff.
\newblock Formules de cubatures m\'ecaniques \`a coefficients non n\'egatifs.
\newblock {\em Bull. Sci. Math. (2)}, 81:123--134, 1957.

\bibitem{SDPT3}
Kim-Chuan Toh, Michael~J. Todd, and Reha~H. T\"ut\"unc\"u.
\newblock On the implementation and usage of {SDPT}3---a {M}atlab software
  package for semidefinite-quadratic-linear programming, version 4.0.
\newblock In {\em Handbook on semidefinite, conic and polynomial optimization},
  volume 166 of {\em Internat. Ser. Oper. Res. Management Sci.}, pages
  715--754. Springer, New York, 2012.

\bibitem{WK06}
X.~Wan and G.~E. Karniadakis.
\newblock Long-term behavior of polynomial chaos in stochastic flow
  simulations.
\newblock {\em Comput. Methods Appl. Mech. Engrg.}, 195(41-43):5582--5596,
  2006.

\bibitem{WW95}
Grzegorz~W. Wasilkowski and Henryk Wo\'zniakowski.
\newblock Explicit cost bounds of algorithms for multivariate tensor product
  problems.
\newblock {\em J. Complexity}, 11(1):1--56, 1995.

\bibitem{XH05}
D.~Xiu and J.~S. Hesthaven.
\newblock High-order collocation methods for differential equations with random
  inputs.
\newblock {\em SIAM J. Sci. Comput.}, 27(3):1118--1139, 2005.

\bibitem{XK02}
D.~Xiu and G.~E. Karniadakis.
\newblock The {W}iener-{A}skey polynomial chaos for stochastic differential
  equations.
\newblock {\em SIAM J. Sci. Comput.}, 24(2):619--644, 2002.

\bibitem{Xiu07}
Dongbin Xiu.
\newblock Efficient collocational approach for parametric uncertainty analysis.
\newblock {\em Commun. Comput. Phys.}, 2(2):293--309, 2007.

\bibitem{Xiu09}
Dongbin Xiu.
\newblock Fast numerical methods for stochastic computations: a review.
\newblock {\em Commun. Comput. Phys.}, 5(2-4):242--272, 2009.

\bibitem{ZRTK2012}
Z.~Zhang, B.~Rozovskii, M.~V. Tretyakov, and G.~E. Karniadakis.
\newblock A multistage {W}iener chaos expansion method for stochastic
  advection-diffusion-reaction equations.
\newblock {\em SIAM J. Sci. Comput.}, 34(2):A914--A936, 2012.

\bibitem{ZTRK14}
Z.~Zhang, M.~V. Tretyakov, B.~Rozovskii, and G.~E. Karniadakis.
\newblock A recursive sparse grid collocation method for differential equations
  with white noise.
\newblock {\em SIAM J. Sci. Comput.}, 36(4):A1652--A1677, 2014.

\bibitem{ZTRK15}
Zhongqiang Zhang, Michael~V. Tretyakov, Boris Rozovskii, and George~E.
  Karniadakis.
\newblock Wiener chaos versus stochastic collocation methods for linear
  advection-diffusion-reaction equations with multiplicative white noise.
\newblock {\em SIAM J. Numer. Anal.}, 53(1):153--183, 2015.

\end{thebibliography}

\end{document}